\documentclass[final]{amsart}

\usepackage{amssymb}
\usepackage[foot]{amsaddr}
\usepackage[utf8]{inputenc}
\usepackage[T1]{fontenc}
\usepackage[english]{babel}
\usepackage{xargs}
\usepackage{graphicx, color}
\usepackage[figurewithin=none]{caption,subcaption} \graphicspath{{./}} \usepackage{enumitem}

\usepackage[notcite,notref]{showkeys} \usepackage{hyperref}
\hypersetup{
	hidelinks,
  colorlinks  = true,    urlcolor    = blue,    linkcolor   = blue,    citecolor   = blue,      pdfauthor   = {A. Montero, D. Pellicer, M. Toledo},pdfsubject  = {Chiral Polytopes},pdftitle    = {Chiral extensions of regular toroids},  }

\usepackage[shadow, colorinlistoftodos,textsize=tiny,obeyFinal]{todonotes}
\setuptodonotes{fancyline, backgroundcolor=gray!10,bordercolor=blue}

\newcommandx{\inline}[2][1=]{\todo[inline, #1]{#2}}
\newcounter{hw}
\newcommandx{\homework}[2][1=]{\todo[inline, caption={Homework \thehw} #1]{\stepcounter{hw} #2}}

 \newcommandx{\danielLong}[2][1=Long todo, usedefault]{\todo[inline, color=green!25, caption={\textbf{Daniel:} #1}]{\textbf{Daniel: #1}. #2}}
 \newcommandx{\micaelLong}[2][1=Long todo, usedefault]{\todo[inline, color=red!25,caption={\textbf{Micael:} #1}]{\textbf{Micael: #1}. #2}}
 \newcommandx{\teroLong}[2][1=Long todo, usedefault]{\todo[inline,color=blue!25, caption={\textbf{Tero:} #1}]{\textbf{Tero: #1}. #2} }

\makeatletter
    \providecommand\@dotsep{5}
\makeatother

\usepackage[pagewise, modulo,mathlines]{lineno}

\theoremstyle{plain}
\newtheorem{thm}{Theorem}[section]

\newtheorem{lem}[thm]{Lemma}
\newtheorem{lemma}[thm]{Lemma}
\newtheorem{prop}[thm]{Proposition}

\newtheorem{coro}[thm]{Corollary}

\theoremstyle{definition}
\newtheorem{defn}[thm]{Definition}

\theoremstyle{remark}

\newtheorem{rem}[thm]{Remark}
\newtheorem{remark}[thm]{Remark}

\numberwithin{equation}{section}

\usepackage{xargs}

\renewcommand{\leq}{\leqslant} \renewcommand{\geq}{\geqslant}
\renewcommand{\epsilon}{\varepsilon} \renewcommand{\subset}{\subseteq}  
\renewcommand{\{}{\lbrace}
\renewcommand{\}}{\rbrace}
\newcommand{\sm}{\setminus} 
\renewcommand{\bar}{\overline}
\renewcommand{\hat}{\widehat}

\newcommand{\bE}{\mathbb{E}}
\newcommand{\bN}{\mathbb{N}}

\newcommand{\bZ}{\mathbb{Z}}

\newcommand{\cF}{\mathcal{F}}

\newcommand{\cK}{\mathcal{K}}
\newcommand{\cM}{\mathcal{M}}

\newcommand{\cP}{\mathcal{P}}

\newcommand{\cS}{\mathcal{S}}
\newcommand{\cT}{\mathcal{T}}
\newcommand{\cU}{\mathcal{U}}

\newcommand{\LL}{\Lambda}
\newcommand{\bLL}{\mathbf{\Lambda}}

\newcommand{\cyvec}[1]{{\mathrm{#1}}}
\newcommand{\vx}{\cyvec{x}}
\newcommand{\vy}{\cyvec{y}}

\newcommand{\va}{\cyvec{a}}

\newcommandx{\rr}[1][1=\ell]{\rho^{#1}}
\newcommand{\id}{\varepsilon}

\DeclareMathOperator{\norm}{Norm}
\DeclareMathOperator{\sgn}{sgn}
\DeclareMathOperator{\aut}{Aut} \DeclareMathOperator{\autp}{\aut^{+}}

\DeclareMathOperator{\rk}{rk}

\DeclareMathOperator{\GPR}{\mathcal{G}}
\DeclareMathOperator{\fl}{\mathcal{F}}
\DeclareMathOperator{\Fw}{\mathcal{F}^{w}}
\DeclareMathOperator{\mon}{Mon}
\DeclareMathOperator{\monp}{\mon^{+}}
\DeclareMathOperator{\monw}{\mon^{w}}

\newcommand{\vect}[1]{\bar{\mathrm{#1}}}

\newcommandx{\dtwoSM}[2][1=\cM, 2=s] {\hat{2}#2^{#1 - 1}}

\newcommandx{\gpr}[2][1=\cK, 2=s]{\GPR_{#2}(#1)}
\newcommand{\baseFlag}{\Phi_{0}}
\newcommandx{\Proot}[2][1=\cP, 2= \baseFlag, usedefault]{\left( #1,#2 \right)}

\newcommand{\cyctwo}[1][n]{ C_{2}^{#1} }

\newcommandx{\groupExt}[1][1=p]{\Gamma^{#1}}

\DeclareMathOperator{\tras}{T}

\newcommand{\kng}{{h}}
\newcommand{\etab}{\bar{\kng}}
\newcommand{\mub}{\bar{\mu}}

\newcommandx{\toroidI}[2][1=n-2, 2=a]{\{4,3^{#1},4\}_{(#2, 0, \dots, 0)}}
\newcommandx{\toroidII}[2][1=n-2, 2=a]{\{4,3^{#1},4\}_{(#2, #2, 0,\dots, 0)}}
\newcommandx{\toroidIII}[2][1=n-2, 2=a]{\{4,3^{#1},4\}_{(#2, #2, \dots, #2)}}
\newcommand{\typetor}[1][n-2]{\{4,3^{#1},4\}}
\newcommand{\E}[1][n]{\bE^{#1}}
\newcommand{\UoverL}{\cU/\bLL}
\newcommand{\minid}{\chi}
\newcommandx{\cl}[1][1=1]{\Lambda_{(#1, 0, \dots, 0)}}
\newcommandx{\fcl}[1][1=1]{\Lambda_{(#1, #1, 0, \dots, 0)}}
\newcommandx{\bcl}[1][1=1]{\Lambda_{(#1, \dots, #1)}}
\newcommandx{\Bcl}[1][1=1]{\mathbf{\cl[#1]}}
\newcommandx{\Bfcl}[1][1=1]{\mathbf{\fcl[#1]}}
\newcommandx{\Bbcl}[1][1=1]{\mathbf{\bcl[#1]}}

\newcommand{\te}{\xi}

\newcommand{\es}{\varsigma}

 \usepackage[style=numeric, sorting=nyt,  maxbibnames=9]{biblatex}
\begin{document}

\title{Chiral extensions of regular toroids}

\author{Antonio Montero}
\email{antonio.montero@fmf.uni-lj.si}
\address[a]{Faculty of Mathematics and Physics, University of Ljubljana, SI-1000 Ljubljana, Slovenia}

\author{Micael Toledo} 
\address[b]{Institute of Mathematics, Physics and Mechanics, Jadranska 19, SI-1000
Ljubljana, Slovenia}
\email{micaelalexitoledo@gmail.com}

\begin{abstract}
    Abstract polytopes are combinatorial objects that generalise geometric objects such as convex polytopes, maps on surfaces and tiling of space. Chiral polytopes are those abstract polytopes that admit full combinatorial rotational symmetry but do not admit reflections. 
    In this paper we build chiral polytopes whose facets (maximal faces) are isomorphic to a prescribed regular cubic tessellation of the $n$-dimensional torus ($n \geq 2$). 
    As a consequence, we prove that for every $d \geq 3$ there exist infinitely many chiral $d$-polytopes.
\end{abstract}

\maketitle 

\listoftodos\relax

\section{Introduction}

An abstract polytope $\cP$ is a partially ordered set that generalises the incidence structure of geometric convex polyhedra to higher dimensions. The rank of $\cP$ is the combinatorial equivalent to the geometric notion of dimension.
Abstract polytopes inherit a natural recursive structure from their convex and geometric analogues: just as a cube can be thought of as a family of six squares (objects of dimension $2$) glued together along their edges (of dimension $1$), an abstract polytope $\cP$ of rank $n$ (called an $n$-polytope) can be thought of as a family of $(n-1)$-polytopes glued along 
 their faces of rank $n-2$. 
Those $(n-1)$-polytopes are the \emph{facets} of $\cP$, and whenever all the facets are {isomorphic} to a fixed polytope $\cK$ we say that $\cP$ is an extension of $\cK$.

The problem of determining whether or not a fixed polytope $\cK$ admits an extension has been part of the theory's development since its beginning. In fact in \cite{DanzerSchulte_1982_RegulareInzidenzkomplexe.I} Danzer and Schulte attack this problem for regular polytopes. 
They prove that every non-degenerate regular polytope $\cK$ admits an extension, and this extension is finite if and only if $\cK$ is finite. In \cite{Danzer_1984_RegularIncidenceComplexes}, Danzer proves that every non-degenerate polytope $\cK$, regardless of their symmetry properties, admits an extension that is finite (resp. regular) if and only if $\cK$ is finite (resp. regular). 
In \cite{Schulte_1983_ArrangingRegularIncidence}, Schulte builds a universal regular extension $\cU$ for every regular polytope $\cK$. 
This extension is universal in the sense that any other regular extension of $\cK$ is a quotient of $\cU$. 
In \cite{Pellicer_2009_ExtensionsRegularPolytopes,Pellicer_2010_ExtensionsDuallyBipartite} Pellicer develops several constructions that have as a consequence that every regular polytope admits a regular extension with prescribed conditions on its local combinatorics. 
In particular, this proves that every regular polytope admits an infinite number of non-isomorphic regular extensions.

Chiral abstract polytopes are those that admit full symmetry by (combinatorial) rotations but do not admit reflections.
They were introduced by Schulte and Weiss in \cite{SchulteWeiss_1991_ChiralPolytopes} as a generalisation of Coxeter's twisted honeycombs in \cite{Coxeter_1970_TwistedHoneycombs}.
The notion of a chiral polytope arises naturally from a geometric idea and examples of chiral $3$-polytopes are abundant. 
Many of them have been part of the literature in the context of maps on surfaces.
In \cite{CoxeterMoser_1972_GeneratorsRelationsDiscrete}, it is proved that there are infinitely many chiral maps on the torus. 
On the other hand, it is known that there are no chiral maps on orientable surfaces of genus $2$, $3$, $4$, $5$ and $6$ (see \cite{ConderDobcsanyi_2001_DeterminationAllRegular, Garbe_1969_UberDieRegularen}, for example). 
The smallest non-toroidal chiral map was constructed by Wilson in \cite{Wilson_1978_SmallestNontoroidalChiral}. 
The results obtained by Sherk in \cite{Sherk_1962_FamilyRegularMaps} imply that infinitely many orientable surfaces admit chiral maps. 

Finding examples of chiral polytopes of higher ranks has proved to be a rather difficult problem. 
Some examples of chiral $4$-polytopes have been built from hyperbolic tilings in \cite{ColbournWeiss_1984_CensusRegular$3$, NostrandSchulte_1995_ChiralPolytopesHyperbolic, SchulteWeiss_1994_ChiralityProjectiveLinear}. 
In \cite{SchulteWeiss_1995_FreeExtensionsChiral} Schulte and Weiss developed a technique to build chiral extensions of polytopes, which introduced the first examples of chiral $5$-polytopes. 
However, this technique cannot be applied twice directly, so it cannot be used to build $6$-polytopes. 
Moreover, the construction introduced by Schulte and Weiss gives a locally infinite polytope. 
The first examples of finite chiral $5$-polytopes were constructed by Conder, Hubard and Pisanski in \cite{ConderHubardPisanski_2008_ConstructionsChiralPolytopes} with the use of computational tools. 
In \cite{DAzevedoJonesSchulte_2011_ConstructionsChiralPolytopes} Breda, Jones and Schulte develop a technique to build new finite chiral $n$-polytopes from known finite chiral $n$-polytopes. 
This technique allowed the construction of concrete examples of chiral polytopes of ranks $3$, $4$ and $5$.

It was not until 2010 that Pellicer proved in \cite{Pellicer_2010_ConstructionHigherRank} the existence of chiral polytopes of all ranks higher than $3$. 
His construction is based on finding a chiral extension of a particular regular polytope and can be applied recursively to the minimal regular cover of the resulting chiral extension. 
Unfortunately, because of the nature of this construction, the size of the polytopes obtained grows so fast that they quickly become conceptually intractable, to the point where determining many of their basic structural properties, such as the kind of facets they have, becomes practically impossible at high ranks.

One of the limitations of any construction of chiral extensions of polytopes is that, unlike the constructions for regular extensions, they cannot be applied recursively. 
If $\cP$ is a chiral polytope, its facets can be either chiral or regular, but the $(n-2)$-faces of $\cP$ must be regular (see \cite{SchulteWeiss_1991_ChiralPolytopes}). 
This implies that if $\cP$ is a chiral extension of $\cK$, then $\cK$ is either regular or chiral with regular facets. 

If $\cK$ is chiral with regular facets, a universal chiral extension of $\cK$ \cite{SchulteWeiss_1995_FreeExtensionsChiral} exists. 
In \cite{CunninghamPellicer_2014_ChiralExtensionsChiral}, Cunningham and Pellicer proved that any finite chiral polytope with regular facets admits a finite chiral extension. 
There are examples of orientably regular polytopes that do not admit a chiral extension (see \cite{Cunningham_2017_NonFlatRegular}), but these examples are, in some sense, degenerate. 

In \cite{ConderHubardOReillyRegueiroPellicer_2015_ConstructionChiral4} the authors build chiral $4$-polytopes with symmetric and alternating groups as automorphism groups. 
In \cite{ConderHubardOReillyRegueiro_ConstructionChiralPolytopes_preprint} Conder, Hubard and O'Reilly-Regeiro extend the techniques in \cite{ConderHubardOReillyRegueiroPellicer_2015_ConstructionChiral4} to build chiral polytopes rank higher than $4$ whose automorphism group is alternating or symmetric. 
These polytopes arise naturally as chiral extensions of the simplex.
To the best of our knowledge, this is the only known construction of chiral extensions of regular polytopes. 

An $(n+1)$-toroid is a quotient of a regular tiling of the Euclidian space $\E$ by a lattice group. 
Apart from a few exceptions, every $(n+1)$-toroid possesses the natural structure of an abstract polytope. 
Moreover, regular toroids are well understood: all of them are quotients of a regular tessellation of $\E$; in particular, if $n\not \in  \{2,4\}$, all the regular $(n+1)$-toroids are a quotient of a cubic tessellation.
The family of cubic regular toroids is arguably the most natural infinite family of regular polytopes of any given rank.
In this paper, we attack the problem of building chiral extensions of cubic regular toroids. 
More precisely, we prove the following result (see Theorem\nobreakspace \ref {thm:elChido} for a more detailed version).

\begin{thm}\label{thm:elNoTanChido}
Let $n \geq 2$. For all but finitely many cubic regular $(n+1)$-toroids $\cT$ there exists a chiral $(n+2)$-polytope whose facets are isomorphic to $\cT$. 
\end{thm} \section{Preliminaries}
\label{sec:basics}

\subsection{Regular abstract polytopes.}\label{sec:HSAP}

Abstract polytopes are structures derived from combinatorial properties of geometric polytopes. 
They generalise convex polytopes, tilings of the Euclidean spaces, maps on surfaces, among others.
Formally, an \emph{abstract polytope of rank $n$} or an \emph{$n$-polytope}, for short, is a partially ordered set $(\cP, \leq)$ (we usually omit the order symbol) that satisfies the properties in Items\nobreakspace  \ref {item:maxmin} to\nobreakspace  \ref {item:diamond}  below. 
 These properties are precisely P1, P2, P3' and P4 in \cite[Section 2A]{McMullenSchulte_2002_AbstractRegularPolytopes}, where they are described in detail.
 \begin{enumerate}
  \item \label{item:maxmin} $\cP$ has a minimum element $F_{-1}$ and a maximum element $F_{n}$.
\item \label{item:flags} Every flag (maximal chain) of $\cP$ contains exactly $n+2$ elements, including $F_{-1}$ and $F_{n}$.
\item \label{item:sfc} $\cP$ is \emph{strongly flag connected}.
 \item \label{item:diamond} $\cP$ satisfies the \emph{diamond condition}.
\end{enumerate}
The elements of $\cP$ are called \emph{faces}. 
We say that two faces $F$ an $G$ are \emph{incident} if $F \leq G$ or $G\leq F$.
The condition in Item\nobreakspace \ref {item:flags} allows us to define a \emph{rank function} $\rk: \cP \to \left\{ -1, \dots, n \right\} $ by $\rk(F)=i$, where $i+2$ is the length of a maximal chain of $\cP$ whose largest face is $F$. 
In particular $\rk(F_{-1}) = -1$ and $\rk(F_{n})=n$. 
We usually call \emph{vertices}, \emph{edges} and \emph{facets} the elements of rank $0$, $1$ and $n-1$ respectively.
In general, a face of rank $i$ is called an $i$-face.
The diamond condition implies that given $i \in \left\{ 0, \dots, n-1 \right\} $, for every flag $\Phi$ of $\cP$ there exists a unique flag $\Phi^{i}$ such that $\Phi$ and $\Phi^{i}$ differ exactly in the face of rank $i$. 
In this situation we say that $\Phi$ and $\Phi^{i}$ are \emph{adjacent} or \emph{$i$-adjacent}, if we need to emphasise $i$.
If $i_{1}, i_{2}, \dots, i_{k} \in \left\{ 0, \dots, n-1 \right\} $, we define recursively $\Phi^{i_{1}, \dots, i_{k}} = (\Phi^{i_{1}, \dots, i_{k-1}})^{i_{k}}$.
This notion of adjacency turns the set $\cF(\cP)$ of flags of $\cP$ into a properly $n$-edge-coloured graph, called the \emph{flag graph of $\cP$}.

For every $i \in \{0,\ldots,n-1\}$, let $r_i$ be the permutation of the flags of $\cP$ mapping every flag $\Phi$ to its $i$-adjacent flag $\Phi^i$. We call the group $\mon(\cP) := \langle r_0, \ldots, r_{n-1} \rangle$  the \emph{monodromy group} of $\cP$. We will let $\mon(\cP)$ act of the left, so that $r_i\Phi$ denotes the image of $\Phi$ under $r_i$. Note that the action of $\mon(\cP)$ is transitive on the set of flags of $\cP$ but it may not be free.

The {\em even monodromy group} of $\cP$ is the subgroup of  $\mon(\cP)$  defined as $\mon^+(\cP) = \langle s_1, \ldots, s_{n-1} \rangle$ where $s_{i}=r_{i-1}r_{i}$ for all $i \in \{1\ldots,n-1\}$. Note that $\mon^+(\cP)$ is a subgroup of index at most $2$ in $\mon(\cP)$. This group will play an important role throughout this paper.

If $F$ and $G$ are two faces of a polytope such that $F \leq G$, the \emph{section} $G/F$ is the restriction of the order of $\cP$ to the set $\left\{ H \in \cP : F \leq H \leq G \right\} $.
Note that if $\rk(F) = i$ and $\rk(G) = j$, the section $G/F$ is an abstract polytope of rank $j-i-1$.
If $F_{0}$ is a vertex, the \emph{vertex-figure} at $F_{0}$ is the section $F_{n}/F_{0}$. 
We sometimes identify each face $F$ with the section $F/F_{-1}$. 
In particular, every facet $F_{n-1}$ can be identified with the section $F_{n-1}/F_{-1}$ of rank $n-1$.

Given an abstract $n$-polytope $\cK$, an $(n+1)$-polytope $\cP$ is an \emph{extension} of $\cK$ if all the facets of $\cP$ are isomorphic to $\cK$.

For $i \in \left\{ 1, \dots, n-1 \right\} $, if $F$ is an $(i-2)$-face and $G$ is an $(i+1)$-face with $F \leq G$ then the section $G/F$ is a $2$-polytope. 
Therefore $G/F$ is isomorphic to a $p_{i}$-gon for some $p_{i} \in \left\{ 2, \dots, \infty \right\} $. 
If the number $p_{i}$ does not depend on the particular choice of $F$ and $G$ but only on $i$ we say that $\cP$ has \emph{Schläfli symbol} $\left\{ p_{1}, \dots, p_{n-1} \right\}$;
in this situation sometimes we just say that $\cP$ is of \emph{type} $\left\{ p_{1}, \dots, p_{n-1} \right\} $.
Note that if $\cP$ is an $n$-polytope of type $\left\{ p_{1}, \dots, p_{n-1} \right\} $, then the the facets of $\cP$ are of type $\left\{ p_{1}, \dots, p_{n-2} \right\} $.
In particular, if $\cP$ is an extension of $\cK$, and $\cP$ has a well-defined Schläfli symbol, all but the last entry of this symbol are determined by $\cK$.

An \emph{automorphism} of an abstract polytope $\cP$ is an order-preserving bijection $\gamma: \cP \to \cP$.
The group of automorphisms of $\cP$ is denoted by $\aut(\cP)$.
The group $\aut(\cP)$ acts naturally on $\fl(\cP)$. We will consider this to be a right action so that \[(r_i\Psi) \gamma = r_i( \Psi \gamma )\] for every flag $\Psi$, $r_i \in \mon(\cP)$ and $\gamma \in \aut(\cP)$.
That is, we may think of $\aut(\cP)$ as the group of all flag permutations that commute with the action of $\mon(\cP)$. As a consequence of the strong-flag-connectivity, the action  of $\aut(\cP)$ on $\fl(\cP)$ is free, although it may not be transitive.

Let $\baseFlag = \left\{ F_{-1}, \dots, F_{n} \right\} $ be a \emph{base flag} of $\cP$  such that $\rk(F_{i}) = i$.
Let $\Gamma \leq \aut(\cP)$ and for $I \subset \left\{ 0, \dots, n-1 \right\} $ let $\Gamma_{I}$ denote the set-wise stabiliser of the chain $\left\{ F_{i} : i \not\in I \right\} \subset \baseFlag $. 
Note that for every pair of subsets $I, J \subset \left\{ 0, \dots, n-1 \right\} $ we have 
\begin{equation}\label{eq:intProperty}
	\Gamma_{I} \cap \Gamma_{J} = \Gamma_{I \cap J}
\end{equation}
We call this condition the \emph{intersection property} for $\Gamma$.

An abstract polytope is $\emph{regular}$ if the action of $\aut(\cP)$ on $\fl(\cP)$ is transitive (hence, regular).
Traditionally, among abstract polytopes the regular ones have been studied most frequently. 
Extensive theory can be found in \cite{McMullenSchulte_2002_AbstractRegularPolytopes}.

A \emph{rooted polytope} is a pair $(\cP, \baseFlag)$ where $\cP$ is a polytope and $\baseFlag$ is a fixed base flag.
If $\cP$ is regular, every two flags are equivalent under $\aut(\cP)$ and the choice of a particular base flag plays no relevant role. 
However, if $\cP$ is not regular, then the choice of the base flag is important.
See \cite{CunninghamPellicer_2018_OpenProblems$k$} for a discussion on rooted $k$-orbit polytopes.
 
If $\Proot$ is a regular rooted polytope, then for every $i \in \left\{ 0, \dots, n-1 \right\} $ there exists an automorphism $\rho_{i}$ such that \[\baseFlag\rho_{i} = \baseFlag^{i}.\]
We call the automorphisms $\rho_{0}, \dots, \rho_{n-1}$ the \emph{abstract reflections} (with respect to the base flag $\baseFlag$). 
It is easy to see that if $\cP$ is a regular $n$-polytope, then $\aut(\cP) = \left\langle \rho_{0}, \dots, \rho_{n-1} \right\rangle $.
It is important to remark that the group elements depend on $\baseFlag$. 
However, since $\aut(\cP)$ is transitive on flags, a change in the base flag corresponds to applying an inner group-automorphism to the generators of $\aut(\cP)$.
More precisely, let $\Phi$ and $\Psi$ be flags of a regular $n$-polytope $\cP$ and let $\rho_{0}, \dots, \rho_{n-1}$ and $\rho'_{0}, \dots, \rho'_{n-1}$ denote the abstract reflections with respect to $\Phi$ and $\Psi$ respectively. 
If $\gamma \in \aut(\cP)$ is such that $\Phi \gamma  = \Psi$, then $\rho'_{i} = \gamma^{-1} \rho_{i} \gamma$.

Note that every regular polytope has a well-defined Schläfli symbol.
If $\cP$  is a regular $n$-polytope of type $\left\{ p_{1}, \dots, p_{n-1} \right\} $, then the abstract reflections satisfy 
\begin{equation}\label{eq:relsRhos}
	\begin{aligned} 
		\rho_{i}^{2} &= \id && \text{for } i \in \left\{ 0, \dots, n-1 \right\} ,\\
		\left( \rho_{i} \rho_{j} \right)^{2} &= \id && \text{if } |i-j| \geq 2, \\
		\left( \rho_{i-1} \rho_{i}\right)^{p_{i}} &= \id && \text{for } i \in \left\{ 1, \dots, n-1 \right\}. 
	\end{aligned}
\end{equation}

If $\baseFlag = \left\{ F_{-1}, \dots, F_{n} \right\} $, the stabiliser of the chain $\left\{ F_{i} : i \not \in I \right\} $ is the group $\left\langle \rho_{i} : i \in I \right\rangle $.
It follows that for regular polytopes the intersection property in Equation\nobreakspace \textup {(\ref {eq:intProperty})} for $\aut(\cP)$ itself is equivalent to
\begin{equation}\label{eq:interPropReg}
 \left\langle \rho_{i} : i \in I \right\rangle \cap \left\langle \rho_{j}: j \in J \right\rangle = \left\langle \rho_{k} : k \in I \cap J \right\rangle   
\end{equation}
for every pair of sets $I, J \subset \left\{ 0, \dots, n-1 \right\} $.

A \emph{string C-group} is a group $\left\langle \rho_{0}, \dots, \rho_{n-1} \right\rangle $ satisfying Equations\nobreakspace \textup {(\ref {eq:relsRhos})} and\nobreakspace  \textup {(\ref {eq:interPropReg})}. 
Clearly, the automorphism group of an abstract polytope is a string C-group.
One of the most remarkable facts in the theory of highly symmetric polytopes is the correspondence between string C-groups and abstract regular polytopes. 
To be precise, for every string C-group $\Gamma$, there exists an abstract regular polytope $\cP = \cP(\Gamma)$ such that $\aut(\cP) = \Gamma$.

The abovementioned correspondence has been used to build families of abstract regular polytopes with prescribed properties. 
For instance, some universal constructions are explored in \cite{Schulte_1983_ArrangingRegularIncidence} and \cite{Schulte_1985_ExtensionsRegularComplexes}. 
In another direction, in \cite{CameronFernandesLeemansMixer_2017_HighestRankPolytope, FernandesLeemans_2018_CGroupsHigh, LeemansMoerenhoutOReillyRegueiro_2017_ProjectiveLinearGroups, Pellicer_2008_CprGraphsRegular} some abstract regular polytopes with prescribed (interesting) groups are investigated. 
Of particular interest for this paper is the work in \cite{Danzer_1984_RegularIncidenceComplexes, Pellicer_2009_ExtensionsRegularPolytopes, Pellicer_2010_ExtensionsDuallyBipartite}, where the problem of finding regular extensions of regular polytopes is addressed. 

\subsection{Rotary polytopes} \label{sec:rotary}

Abstract regular polytopes are those with a maximal degree of reflectional symmetry. 
A slightly weaker symmetry condition than regularity for an abstract polytope is to admit all possible rotational symmetries.
In a similar way as has been done for maps (see \cite{Wilson_1978_RiemannSurfacesOver}, for example), we call these polytopes \emph{rotary polytopes}. 
In this section, we review some of the theory of rotary polytopes.
Most of this theory is developed in \cite{SchulteWeiss_1991_ChiralPolytopes}.

An abstract polytope $\cP$ is \emph{orientable} if $\mon^+(\cP) = \langle r_{i-1}r_{i} : 1 \leq i \leq n-1 \rangle$ has index $2$ in $\mon(\cP)$;
otherwise, $\cP$ is \emph{non-orientable}.
If $\Proot$ is a rooted orientable polytope, we may bicolour the flags of $\cP$ by defining $\baseFlag$ as white and recursively colouring any other flag $\Psi$  black (resp. white) if and only if it is adjacent to a white (resp. black) flag.
We denote the set of white flags by $\Fw(\cP)$.
This is just another way of naming what in \cite{SchulteWeiss_1991_ChiralPolytopes} are called \emph{even} and \emph{odd} flags. In fact, the set of white flags corresponds to the orbit of the base $\Phi_0$ under the action of the even mondoromy group $\mon^+(\cP)$.
For convenience, if $\cP$ is non-orientable, we set  $\Fw(\cP) = \fl (\cP)$.
In other words, every flag is white.

If $\cP$ is an abstract polytope, the \emph{rotational group} $\autp(\cP) $ of $\cP$ is the subgroup of $\aut(\cP)$ that permutes the set of white flags. 
A polytope is \emph{rotary} if $\autp(\cP)$ acts transitively on the set of white flags. 
It is clear that a rotary non-orientable polytope is a regular polytope. 
Therefore, we restrict our discussion below to orientable polytopes.
Note that the choice of the base flag of $\cP$ plays a stronger role now. 
In particular, it defines the set of white flags.
In the discussion below it is assumed that $\cP$ is actually a rooted polytope $\Proot$.

If $\Proot$ is a  polytope, for every $i \in \left\{ 1, \dots, n-1 \right\} $ the flag $\baseFlag^{i,i-1}$ is a white flag.
Therefore, if $\cP$ is a rotary polytope, there exists an automorphism $\sigma_{i}$ such that \[\baseFlag \sigma_{i} = \baseFlag^{i, i-1}.\]
The automorphisms $\sigma_{1}, \dots, \sigma_{n-1}$ are called the \emph{abstract rotations} with respect to $\baseFlag$.
It is easy to see that $\autp(\cP) = \left\langle \sigma_{1}, \dots, \sigma_{n-1} \right\rangle $. 
We emphasise that the abstract rotations depend on the choice of the base flag.

If $\cP$ is a rotary $n$-polytope, then $\cP$ has a well-defined Schläfli symbol. 
If $\cP$ is of type $\left\{ p_{1}, \dots, p_{n-1} \right\} $ the automorphisms $\sigma_{1}, \dots, \sigma_{n-1}$ satisfy 
\begin{equation}\label{eq:relsSigmas}
 \begin{aligned}
		 	  \sigma_{i}^{p_{i}} &= \id, && \text{and} \\
	  (\sigma_{i} \sigma_{i+1} \cdots \sigma_{j})^{2} &= \id && \text{for}\ 1 \leq i < j \leq n-1.
 \end{aligned}
\end{equation}

Sometimes it is useful to consider an alternative set of generators for $\autp(\cP)$. 
For $i,j \in \{0, \dots, n-1\}$ with $i <j$ we define the automorphisms 
\[\tau_{i,j} = \sigma_{i+1}\cdots \sigma_{j}.\] 
Note that this is a small change with respect to the notation of \cite[Eq. 5]{SchulteWeiss_1991_ChiralPolytopes}. 
What they call $\tau_{i,j}$ for us is $\tau_{i-1,j}$. 
Observe that $\tau_{i-1,i} = \sigma_{i}$ for $i \in \{1, \dots, n-1\}$. It is also convenient to define $\tau_{j,i} = \tau_{i,j}^{-1}$ for $i < j$ and $\tau_{-1,j}= \tau_{i,n} = \tau_{i,i} = \id$ for every $i,j \in \{0, \dots, n-1\}$. 
In particular, we have that $\langle \tau_{i,j} : i,j \in \{0, \dots, n-1\} \rangle = \langle \sigma_{1}, \dots, \sigma_{n-1} \rangle$. 
We also have \[\baseFlag\tau_{i,j}= \baseFlag^{j,i}.\]

Moreover, if $\baseFlag = \left\{ F_{-1}, \dots, F_{n} \right\} $, and $I \subset \left\{ 0, \dots, n-1 \right\} $,  the stabiliser of the chain $\left\{ F_{i} : i \not\in I  \right\} $ is the group $\left\langle \tau_{i,j} : i,j \in I \right\rangle $.
It follows that the intersection property in Equation\nobreakspace \textup {(\ref {eq:intProperty})} for $\autp(\cP)$ can be written as
\begin{equation}\label{eq:intPropertyChiral}
 \left\langle \tau_{i,j} : i,j \in I \right\rangle \cap \left\langle \tau_{i,j} : i,j \in J \right\rangle = \left\langle \tau_{i,j} : i,j \in I \cap J \right\rangle.  
\end{equation}

If $\cP$ is a regular polytope with automorphism group $\aut(\cP) = \left\langle \rho_{0}, \dots, \rho_{n-1} \right\rangle $, then $\cP$ is rotary with $\autp(\cP) = \left\langle \sigma_{1}, \dots, \sigma_{n-1} \right\rangle $ where $\sigma_{i} = \rho_{i-1} \rho_{i}$. 
If $\cP$ is also orientable, we say that $\cP$ is \emph{orientably regular}.
In this situation $\autp(\cP)$ is a proper subgroup of $\aut(\cP)$ of index $2$.
Furthermore, $\autp(\cP)$ induces two flag-orbits, namely, the white flags and the black flags.

If $\cP$ is rotary but not regular, then $\aut(\cP) = \autp(\cP)$ and this group induces precisely two orbits in flags in such a way that adjacent flags belong to different orbits. 
In this case we say that $\cP$ is \emph{chiral}. 
Chiral polytopes were introduced by Schulte and Weiss in \cite{SchulteWeiss_1991_ChiralPolytopes} as a combinatorial generalisation of Coxeter's twisted honeycombs in \cite{Coxeter_1970_TwistedHoneycombs}.

If $(\cP,\baseFlag)$ is a rooted chiral polytope, the \emph{enantiomorphic form of $\cP$}, denoted by $\bar{\cP}$, is the rooted polytope $(\cP, \baseFlag^{0})$.
In the classic development of the theory of chiral polytopes, the enantiomorphic form of $\cP$ is usually thought as the mirror image of $\cP$, thus as a polytope which is different from but isomorphic to $\cP$.
However, when treated as rooted polytopes, it is clear that the only difference is the choice of the base flag.
The underlying partially ordered set is exactly the same. 
For a traditional but detailed discussion about enantiomorphic forms of chiral polytopes we suggest \cite[Section 3]{SchulteWeiss_1994_ChiralityProjectiveLinear}.

The automorphism group of $\bar{\cP}$ is generated by the automorphisms $\sigma'_{1}, \dots, \sigma'_{n-1}$ where $(\baseFlag^{0}) \sigma'_{i} = (\baseFlag^{0})^{i,i-1}$.
It is easy to verify that $\sigma'_{1} = \sigma_{1}^{-1} $, $\sigma_{2} = \sigma_{1}^{2} \sigma_{2}$ and for $i \geq 3$, $\sigma'_{i} = \sigma_{i}$. 
If $\cP$ is orientably regular the conjugation by $\rho_{0}$ defines a group automorphism $\rho: \autp(\cP) \to \autp(\cP)$ that maps $\sigma_{i}$ to $\sigma'_{i}$. 

A group that satisfies the relations in Equation\nobreakspace \textup {(\ref {eq:relsSigmas})} together with the intersection property in Equation\nobreakspace \textup {(\ref {eq:intPropertyChiral})} must be the rotation group of a rotary polytope.
Moreover, the existence of the group-automorphism $\rho: \autp(\cP) \to \autp(\cP)$ mentioned above determines whether or not the rotary polytope is regular.
More precisely, the following result holds.

\begin{thm}[{\cite[Theorem 1]{SchulteWeiss_1991_ChiralPolytopes}}]\label{thm:chiralGroups}
Let $3\leq n$, $2 < p_{1}, \dots, p_{n-1} \leq \infty$ and let $\Gamma$ be a group such that $\Gamma=\langle \sigma_{1}, \dots, \sigma_{n-1} \rangle$. For every $i,j \in \{-1, \dots, n\}$, with $i \neq j$ define
\[\tau_{i,j}=
\begin{cases}
	\id & \text{if}\ i<j\ \text{and}\ i=-1\ \text{or}\ j=n,\\ 
  \sigma_{i+1} \cdots \sigma_{j} & \text{if}\ 0\leq i<j\leq n-1, \\
  \sigma^{-1}_{j} \cdots \sigma^{-1}_{i+1} & \text{if}\ 0\leq j<i\leq n-1. 
\end{cases}
\]

Assume that $\Gamma$ satisfies the relations in  Equation\nobreakspace \textup {(\ref {eq:relsSigmas})}. Assume also that Equation\nobreakspace \textup {(\ref {eq:intPropertyChiral})} holds. Then
\begin{enumerate}
  \item\label{part:existence} There exists a rotary polytope $\cP=\cP(\Gamma)$ such that $\autp(\cP)=\Gamma$ and $\sigma_{1}, \dots, \sigma_{n-1}$ act as abstract rotations for some flag of $\cP$.
  \item\label{part:facets} $\cP$ is of type $\{p_{1}, \dots, p_{n-1}\}$. The facets and vertex-figures of $\cP$ are isomorphic to $\cP(\langle\sigma_{1}, \dots, \sigma_{n-2}\rangle)$ and $\cP(\langle \sigma_{2}, \dots, \sigma_{n-1})\rangle$, respectively. In general if $n \geq 4$, $F$ is a $(k-2)$-face and $G$ is an incident $(l+1)$-face, for  $1 \leq k < l \leq n-1$,  the section $G/F$ is a rotary $(l-k+2)$-polytope isomorphic to $\cP(\langle \sigma_{k} \dots, \sigma_{l} \rangle)$. 
  \item\label{part:chirality} $\cP$ is orientably regular if and only if there exists an involutory group automorphism $\rho:\Gamma \to \Gamma$ such that $\rho: \sigma_{1} \mapsto \sigma^{-1}_{1}$, $\rho: \sigma_{2} \mapsto \sigma_{1}^{2}\sigma_{2}$ and $\rho: \sigma_{i}\mapsto \sigma_{i}$ for $i\geq 3$.
\end{enumerate}
\end{thm}

\subsection{Extensions of rotary polytopes}
Recall that if $\cP$ is a rotary $n$-polytope, then its facets must be rotary (chiral or orientably regular) but its $(n-2)$-faces must be orientably regular (see \cite[Proposition 9]{SchulteWeiss_1991_ChiralPolytopes}).
It follows that if $\cP$ is a chiral extension of a polytope $\cK$, then $\cK$ is either orientably regular or chiral with regular facets. 
We shall carry these assumptions over $\cK$ through out this section.

A natural approach to build extensions of rotary polytopes is to use Part\nobreakspace \ref {part:facets} of Theorem\nobreakspace \ref {thm:chiralGroups}.
More precisely, let $\cK$ be a rotary $n$-polytope with $\autp(\cK)=\left\langle \sigma_{1}, \dots, \sigma_{n-1} \right\rangle $.
Assume that $\Gamma = \left\langle \bar{\sigma}_{1}, \dots, \bar{\sigma}_{n-1}, \bar{\sigma}_{n} \right\rangle $ is a group satisfying Equation\nobreakspace \textup {(\ref {eq:relsSigmas})} (for rank $n+1$) and such that the mapping $\sigma_{i} \mapsto \bar{\sigma_{i}}$ (for $1 \leq i \leq n-1 $) is an embedding of $\autp(\cK)$ into $\Gamma$.
In order to build a polytope from $\Gamma$ we need to prove that it has the intersection property (Equation\nobreakspace \textup {(\ref {eq:intPropertyChiral})}). 
In this situation the polytope $\cP:=\cP(\Gamma)$ obtained from $\Gamma$ is rotary.
If $\cK$ is chiral, then $\cP$ must be chiral; 
however, if $\cK$ is regular then $\cP(\Gamma)$ might be orientably regular.
To prove that $\cP$ is chiral we need to show that $\Gamma$ does not admit a group automorphism as the one described in Part\nobreakspace \ref {part:chirality} of Theorem\nobreakspace \ref {thm:chiralGroups}.

In this section we develop a series of small results that will prove to be useful in Section\nobreakspace \ref {sec:chirExtRegToroids}, when we establish our main results an build chiral extensions of regular toroids. 
Most of these results are either straightforward or can be found in the literature, hence we omit the proofs and rather provide the appropriate references.
Moreover, instead of assuming that the group $\autp(\cK)$ is embedded into a group $\Gamma$, we shall abuse notation and simply denote by $\sigma_1, \dots, \sigma_n$ the set of distinguished generators of $\Gamma$ and think of $\autp(\cK)$ as the subgroup $ \Gamma_{n}= \langle \sigma_1, \dots, \sigma_{n-1} \rangle$.

We begin with a result that guarantees that the group $\Gamma$ has the appropriate relations.

\begin{lem}\label{lem:relsOfExtensions}
  Let $\Gamma = \langle \sigma_{1}, \dots, \sigma_{n} \rangle$ be a group with the property that the subgroup $\Gamma_{n} = \langle \sigma_{1}, \dots, \sigma_{n-1} \rangle$ satisfies Equation\nobreakspace \textup {(\ref {eq:relsSigmas})}. If the group elements $\sigma_{1}, \dots, \sigma_{n}$ satisfy the equations
\begin{equation}\label{eq:secondSetOfRels}
    \begin{aligned}
      \sigma_{n}^{p_{n}} &= \id, && \text{for some}\ p_{n} \geq 3,\\
      (\sigma_{i} \cdots \sigma_{n})^{2} &= \id, && \text{for}\ 1\leq i \leq n-1,
    \end{aligned}
  \end{equation}
  then $\Gamma$ itself satisfies Equation\nobreakspace \textup {(\ref {eq:relsSigmas})}.
\end{lem}

There is a similar result regarding the intersection property. 
We give it here without proof. 
This is \cite[Lemma 10]{SchulteWeiss_1991_ChiralPolytopes}.

\begin{lem}\label{lem:Lemma10}
  Let $n \geq 3$ and $\Gamma = \langle \sigma_{1}, \dots, \sigma_{n} \rangle$ be a group which satisfies Equation\nobreakspace \textup {(\ref {eq:relsSigmas})}. Assume that the subgroup $\Gamma_{n-1} = \langle \sigma_{1}, \dots, \sigma_{n-1} \rangle$ has the intersection property in Equation\nobreakspace \textup {(\ref {eq:intPropertyChiral})} with respect to its generators. Also, suppose that the following intersection conditions hold:
  \begin{equation}\label{eq:Lemma10}
    \langle \sigma_{1}, \dots, \sigma_{n-1} \rangle \cap \langle \sigma_{j}, \dots, \sigma_{n} \rangle = \langle \sigma_{j}, \dots, \sigma_{n-1} \rangle
  \end{equation}
	for $j \in \{2, \dots, n\}$. Then $\Gamma$ itself has the intersection property of Equation\nobreakspace \textup {(\ref {eq:intPropertyChiral})}.
\end{lem}

Sometimes it is useful to consider a different set of generators for $\Gamma$. 
For $i \in \{1, \dots, n\}$ define
\begin{equation}\label{eq:goodTaus}
  \tau_{i}=\sigma_{1}\cdots \sigma_{i}.
\end{equation}
Assume for a second that $\Gamma = \autp(\cP)$ for a rotary polytope $\cP$ and that $\hat{\Phi_{0}}$ denotes the base flag of $\cP$, then \[\hat{\Phi_{0}} \tau_{i} = \hat{\Phi_{0}}^{i,0}.\] 
We also have $\sigma_{1} = \tau_{1}$ and for $2 \leq i \leq n-1$, $\sigma_{i} = \tau^{-1}_{i-1} \tau_{i}$. 
It follows that 
$\Gamma_{n}=\langle \tau_{1}, \dots, \tau_{n-1} \rangle$.
The following result is essentially the same as Lemma\nobreakspace \ref {lem:relsOfExtensions} but in terms of the generators $\tau_{1}, \dots, \tau_{n}$.

\begin{coro}\label{coro:relsOfExtensionsTaus}
  Let $\Gamma = \langle \sigma_{1}, \dots, \sigma_{n} \rangle$ be a group with the property that the subgroup $\Gamma_{n} = \langle \sigma_{1}, \dots, \sigma_{n-1} \rangle$ satisfies Equation\nobreakspace \textup {(\ref {eq:relsSigmas})}. For $i \in \{1, \dots, n\}$, let $\tau_{i} = \sigma_{1} \cdots \sigma_{i}$. Then the set of relations of Equation\nobreakspace \textup {(\ref {eq:secondSetOfRels})} is equivalent to the set of relations
  \begin{equation}\label{eq:relsOfExtensionsTaus}
    \begin{aligned}
      (\tau_{n-1}^{-1}\tau_{n})^{p_{n}} &= \id,\\
      \tau^{2}_{n} &= \id,\\
      (\tau^{-1}_{i} \tau_{n})^{2} &= \id, && \text{for}\ i \in \{1, \dots, n-2\}.
    \end{aligned}
  \end{equation}
\end{coro}

Finally, in order to determine if the group $\Gamma = \langle \tau_{1}, \dots, \tau_{n} \rangle$ satisfies Corollary\nobreakspace \ref {coro:relsOfExtensionsTaus} and Lemma\nobreakspace \ref {lem:Lemma10} is the automorphism group or a chiral polytope, (and not the rotation subgroup of an orientably regular polytope) we need to show that there is no group automorphism $\alpha:\Gamma \to \Gamma$ such as the one described in Part\nobreakspace \ref {part:chirality} of Theorem\nobreakspace \ref {thm:chiralGroups}. 
The properties of this automorphism are described, in terms of $\tau_{1}, \dots, \tau_{n}$ in the following result.

\begin{coro}\label{cor:gpAutAltGens}
	Let $\Gamma = \langle \sigma_{1}, \dots, \sigma_{n} \rangle$ be a group satisfying Equation\nobreakspace \textup {(\ref {eq:relsSigmas})}. Assume that there is a group automorphism $\alpha: \Gamma \to \Gamma$ satisfying the conditions of Part\nobreakspace \ref {part:chirality} of Theorem\nobreakspace \ref {thm:chiralGroups}, that is, $\alpha: \sigma_{1} \mapsto \sigma_{1}^{-1}$, $\alpha: \sigma_{2} \mapsto \sigma_{1}^{2} \sigma_{2} $ while fixing $\sigma_{i}$ for $i \geq 3$. If $\tau_{i} = \sigma_{1} \cdots \sigma_{i}$, then $\alpha$ satisfies
	\begin{equation}\label{eq:gpAutAltGens}
		\begin{aligned}
			\alpha(\tau_{1}) &= \tau_{1}^{-1}, \\
			\alpha(\tau_{i}) &= \tau_{i} && \text{if}\ i \geq 2.
		\end{aligned}
	\end{equation}
	Conversely, a group automorphism $\alpha:\Gamma \to \Gamma$ satisfying Equation\nobreakspace \textup {(\ref {eq:gpAutAltGens})} also satisfies the conditions of Part\nobreakspace \ref {part:chirality} of Theorem\nobreakspace \ref {thm:chiralGroups}.
\end{coro}

\section{Monodromy group of rotary polytopes}
Let $\cP$ be an abstract polytope and consider its monodromy group $\mon(\cP)$.
It is well-known that if $\cP$ is a regular polytope then $\mon(\cP) \cong \aut(\cP)$ with the isomorphism mapping $r_{i}$ to $\rho_{i}$ (see \cite[Theorem 3.9]{MonsonPellicerWilliams_2014_MixingMonodromyAbstract}). 
If $\cP$ is a rotary polytope, then the action of $\autp(\cP)$ and the action of $\mon(\cP)$ have an interesting relationship.
This relationship is described in the following proposition (see \cite[Lemma 2.5 and Proposition 2.7]{Pellicer_2010_ConstructionHigherRank}).

\begin{prop}\label{prop:autpAsMonPChir}
  Let $\cP$ be a rotary $n$-polytope with base flag $\baseFlag$. Let $\autp(\cP)=\langle \sigma_{1}, \dots,  \sigma_{n-1} \rangle$ and $\mon(\cP)=\langle r_{0}, \dots, r_{n-1} \rangle$ be the rotation and monodromy groups of $\cP$ respectively, and consider the even monodromy group $\monp(\cP) = \langle s_{1}, \dots, s_{n-1} \rangle$ of $\cP$, where $s_{i}=r_{i-1}r_{i}$.
  Then the following hold
\begin{enumerate}
  \item\label{part:jumpsChir} For every $i_{1}, \dots,i_{k} \in \{1, \dots, n-1\}$ \[s_{i_{1}}\cdots s_{i_{k}}\baseFlag=\baseFlag\sigma_{i_{1}}\cdots\sigma_{i_{k}}.\]
  \item\label{part:relationsChir} An element $s_{i_{1}}\cdots s_{i_{k}}\in\monp(\cP)$ fixes the base flag if and only if $\sigma_{i_{1}}\cdots \sigma_{i_{k}} = \epsilon$. 
  In this situation, $s_{i_{1}} \cdots s_{i_{k}}$ stabilises every white flag. 
  \item\label{part:isomorphismChir} If $\Fw(\cP)$ denotes the set of white flags of $\cP$, then there is an isomorphism $f: \monw(\cP) \to \autp(\cP)$ with $\monw(\cP)= \langle \bar{s}_{1}, \dots, \bar{s}_{n-1} \rangle$, where $\bar{s}_{i}$ denotes the permutation of $\Fw(\cP)$ induced by $s_{i}$. 
  This isomorphism maps $\bar{s}_{i}$ to $\sigma_{i}$ for every $i \in \{1, \dots, n-1\}$.
\end{enumerate} 
\end{prop}

Note that unlike the regular case, where there is an isomorphism from $\mon(\cP)$ to $\aut(\cP)$ mapping $r_{i}$ to $\rho_{i}$, if $\cP$ is chiral, the mapping $g: \monp(\cP) \to \autp(\cP)$ defined by $ g: s_{i} \mapsto \sigma_{i}$ is not an isomorphism. 
By Part\nobreakspace \ref {part:relationsChir} of Proposition\nobreakspace \ref {prop:autpAsMonPChir}, this mapping is a well-defined epimorphism, but in general it is not injective. 
In other words, there are non-trivial elements of $\monp(\cP)$ that fix every white flag of $\cP$. 
The subgroup of $\monp(\cP)$ containing all such elements (i.e. the kernel of $g$) is called the \emph{chirality group}. 
For some uses and properties of the chirality group see \cite[Section 3]{Cunningham_2012_MixingChiralPolytopes} and \cite[Section 7]{MonsonPellicerWilliams_2014_MixingMonodromyAbstract}.

To avoid confusion, we will try to avoid the use of $\monp(\cP)$ and instead use the group $\monw(\cP)$, which is the permutation group on $\Fw(\cP)$ induced by the action of $\monp(\cP)$.
Since we will only use the action of $\monp(\cP)$ on white flags,  it is safe to abuse notation and identify $s_{i} = r_{i-1} r_{i} \in \monp(\cP)$ with $\bar{s_{i}} \in \monw(\cP)$, the permutation induced by $s_{i}$ on the set $\Fw(\cP)$.

\section{Geometric cubic $(n+1)$-toroids}
Now we turn our attention to cubic toroids. Most of what is discussed in this section is somewhere in the literature. We shall use \cite[Sections 6A and 6D]{McMullenSchulte_2002_AbstractRegularPolytopes} as our main reference, but \cite{CollinsMontero_2021_EquivelarToroidsFew,HubardOrbanicPellicerWeiss_2012_SymmetriesEquivelar4} are also relevant.
Throughout this section, $\cU$ will denote the cubic honeycomb of the Euclidean space $\bE^{n}$ of type $\typetor$ and $\tras(\cU)$ its translation group. 
A cubic $(n+1)$-toroid is the quotient of $\UoverL$ of $\cU$ by a \emph{lattice group} $\bLL$, that is, a subgroup $\bLL \leq \tras(\cU)$ generated by $n$ linearly independent translations. 
If $\bLL = \left\langle t_{1}, \dots, t_{n} \right\rangle$ and $v_{i}$ is the translation vector associated to $t_{i}$, then  $\{v_{i} : 1 \leq i \leq n \}$ is a \emph{basis} for $\bLL$.
The \emph{lattice} $\LL$ associated with $\bLL$ is the orbit of the origin $o$ of $\E$. That is,
\[\LL= o\bLL = \{m_1 v_1 + \cdots + m_n v_n : m_1, \dots, m_n \in \bZ \}.\]
Clearly the basis $\left\{ v_{1}, \dots, v_{n} \right\}$ determines the lattice  $\LL$ and the lattice group $\bLL$ itself.

Geometrically, we can identify the toroid $\UoverL$ with the corresponding tessellation of the torus $\E/\bLL$.
The (open) \emph{Dirichlet domain} (centered at $o$) associated with $\bLL$ is the set \[ D(\bLL) = \left\{ x \in \E : d(o,x) < d(o t, x) \text{ for all } t \in \bLL\sm \{1\} \right\}. \]
We can geometrically realise $\E/\bLL$ by identifying the points of the closure of $D(\bLL)$ that are equivalent modulo $\bLL$.
Observe that by definition, no two points in $D(\bLL)$ are equivalent under $\bLL$ and such identifications happen in the boundary of $D(\bLL)$.

Now we turn our attention to the symmetries of $\cU$. 
The group of (geometric) symmetries of $\cU$ is generated by the reflections $R_i$, $i\in \left\{ 0, \dots, n \right\}$, on hyperplanes of $\E$ and play the role of the abstract reflections in Equation\nobreakspace \textup {(\ref {eq:relsRhos})}. 
Moreover, the symmetry group $G(\cU)$ coincides with $\aut(\cU)$, which is isomorphic to the affine string Coxeter group $[4, 3^{n-2}, 4]$ (often denoted $\tilde{C}_{n}$). That is, the generating reflections satisfy the defining relations:
\begin{equation}\label{eq:relsRrefCubes}
	\begin{aligned} 
		R_{i}^{2} &= \id && \text{for } i \in \left\{ 0, \dots, n \right\} ,\\
		\left( R_{i} R_{j} \right)^{2} &= \id && \text{if } |i-j| \geq 2, \\
		\left( R_{0} R_{1}\right)^{4} = \left( R_{n-1} R_{n}\right)^{4} &= \id, \\
		\left( R_{i} R_{i+1}\right)^{3} &= \id && \text{for } i \in \left\{ 1, \dots, n-2 \right\}. 
	\end{aligned}
\end{equation}

Up to similarity, we may assume that the vertex set of $\cU$ is the set $(\frac{1}{2}, \dots, \frac{1}{2}) + \bZ^{n}$ so that the center of the facets of $\cU$ are precisely the points of integer coordinates. 
In this situation we may give explicit definitions for the reflections $R_{i}$ (see \cite[Table 2]{CollinsMontero_2021_EquivelarToroidsFew}) so that the group $G_{o}(\cU) = \left\langle R_{0}, \dots, R_{n-1} \right\rangle $ preserves the origin $o$ (acting as the stabiliser of the corresponding facet).
The group $G^{n}_{o}(\cU)=\left\langle R_{1}, \dots, R_{n-1} \right\rangle $ acts as the symmetric group $S_{n}$ permuting the coordinate axes. 
On the other hand, the group $H=\left\langle R_{0}^{G_{o}(\cU)} \right\rangle $ generated by the conjugates of $R_{0}$ under $G_{o}(\cU)$ is the group generated by the $n$ reflections on the coordinate hyperplanes (isomorphic to $C_{2}^{n}$).
In fact, it can be seen that the group $G_{o}(\cU)$ is isomorphic to $G^{n}_{o}(\cU) \ltimes H \cong S_{n} \ltimes C_{2}^{n}$ (with $S_{n}$ acting on $ C_{2}^{n}$ by permutation of coordinates).
The product $T_{1}= R R_{n}$, with $R \in G_{o}(\cU)$ a suitable conjugate of $R_{0}$ such that the hyperplanes of $R_{n}$ and $R$ are parallel, is a translation in $\tras(\cU)$ with respect to the vector $e_{1}$.
The conjugates of $T_{1}$ under $G_{o}(\cU)$ are precisely the translations with respect to the vectors $e_{i}$, $i \in \left\{ 1, \dots, n \right\}$.  
The group generated by those translations acts regularly on the facets of $\cU$, which implies that $\tras(\cU) = \left\langle T_{i} : i \in \left\{ 1, \dots, n \right\} \right\rangle \cong \bZ^{n} $ (with $T_{i}$ the translation with respect to $e_{i}$).
In fact, we have that 
\[G(\cU) \cong \aut(\cU) \cong G_{o}(\cU) \ltimes \tras(\cU) \cong (S_{n}\ltimes C_{2}^{n}) \ltimes \bZ^{n},\]
where for $t \in \bZ^{n}$, an element $\sigma \in S_{n}$ acts on $t$ by permuting coordinates while  $\vx \in  C_{2}^{n}$ acts on $t$ by multiplying by $-1$ the entries of $t$ on the coordinates where $\vx$ is non-trivial.

Now we are interested in which automorphisms of $\cU$ induce automorphisms of $\UoverL$. 
Those automorphisms are discussed in detail in \cite[Sect. 3]{HubardOrbanicPellicerWeiss_2012_SymmetriesEquivelar4} and \cite[Sect. 2]{CollinsMontero_2021_EquivelarToroidsFew}. 
The results obtained there are summarized in the following lemma.

\begin{lemma}\label{lem:toroids}
	With the notation given above, the following statements hold.
	\begin{enumerate}
		\item A symmetry $S \in G(\cU)$ projects to an automorphism of $\UoverL$ if and only if $S$ normalises $\bLL$. 
    This is always the case if $S \in \tras(\cU)$. If $S\in G_{o}(\cU)$, then $S \in \norm(\bLL)$  if and only if $\LL S = \LL$.
		\item Since all lattices are centrally symmetric, $\minid:x \mapsto -x$ always projects to an automorphism of $\UoverL$.
		\item The automorphism group $\aut (\UoverL)$ of $\UoverL$ is isomorphic to \[(K' \ltimes \tras(\cU)	) / \bLL \cong K' \ltimes (\tras(\cU) / \bLL )\] where $K'=\{S \in G_{o}(\cU) : \ S^{-1} \bLL S = \bLL \} =\{S \in G_{o}(\cU) \colon\  \LL S = \LL \}$. 
    In particular $\langle \minid \rangle \leq K' \leq G_{o}(\cU)$. The group $\aut(\UoverL)$ has $k$ orbits on the set of flags of $\UoverL$ if and only if the index of $K'$ in $G_{o}(\cU)$ is $k$.
		\item The toroids $\UoverL$ and $\UoverL'$ are isomorphic if and only if $\bLL$ and $\bLL'$ are conjugate in $G(\cU)$. 
    This in turn is true if and only if  there exists $S \in G_{o}(\cU)$ such that $\LL S = \LL'$.
	\end{enumerate}
\end{lemma}

An immediate consequence of Lemma\nobreakspace \ref {lem:toroids} is that $\UoverL$ is regular if and only if $\LL$ is preserved by $G_{o}(\cU)$.
Those lattices are classified in \cite[Section 6D]{McMullenSchulte_2002_AbstractRegularPolytopes}. Such a lattice must be an integer multiple of one the following:
\begin{itemize}
  \item $\cl$: the lattice associated with $\bZ^{n}$ with basis $\left\{ e_{1}, \dots, e_{n} \right\}$.
  \item $\fcl$: the index-$2$ sublattice of $\cl$ consisting of the points whose coordinate sum is even. The set $\left\{  e_{1} + e_{2}, e_{2}- e_{1}, \dots, e_{n}-e_{n-1}\right\}$ is a basis for $\fcl$. If $n=3$ this lattice is often called the \emph{face-centred cubic lattice}.
  \item $\bcl$:  the index-$2^{n-1}$ sublattice of $\cl$ consisting of the points whose coordinates have the same parity. The set $\left\{ 2e_{1}, \dots, 2e_{n-1}, e_{1}+\cdots+e_{n} \right\}$ is basis for $\bcl$. If $n=3$ this is often called the \emph{body-centred} lattice.
\end{itemize}
Naturally, the corresponding lattice groups are denoted by $\Bcl$, $\Bfcl$ and $\Bbcl$. 
Notice that for $n=2$ we have $\fcl = \bcl$. 
In short, we have the following theorem.

\begin{thm}[ See {\cite[1D,6D1 and 6D4]{McMullenSchulte_2002_AbstractRegularPolytopes}}]
  \label{thm:toroids}
  Let $n \geq 2$ and let $\cT = \UoverL$ be a regular cubic $(n+1)$-toroid. Then $\cT \cong \cU/\mathbf{\LL_{\bar{\va}}}$ where $\bar{\va} = (a^{k}, 0^{n-k})$ for some $a \in \bN$ and $k \in \left\{ 1,2,n \right\}$.
  In every case \[ \aut(\cT) \cong (S_{n} \ltimes C_{2}^{n}) \ltimes \bZ^{n}/\LL_{\bar{\va}}.  \] 
\end{thm} \section{Combinatorial toroids}
\label{sect:combinatorialToroids}

In order to give the appropriate notation for this section we are going to study the structure of the regular toroid $\cT$. We label the flags of $\cT$ as follows. Recall that $\aut(\cT) \cong (S_{n} \ltimes \cyctwo) \ltimes ( \bZ^{n}/\LL_{\bar{\va}})$ for some lattice group $\LL_{\bar{\va}}$ (see Theorem\nobreakspace \ref {thm:toroids}).
Let us denote by $\tras(\cT)$ group $\bZ^{n}/\LL_{\bar{\va}}$ of translations of $\cT$.
Let $\gamma \in \aut(\cT)$ and let $\sigma \in S_{n}$, $\vx \in \cyctwo$ and $t \in \tras (\cT)$ such that $\gamma = \sigma \cdot \vx \cdot t$. 
If $\Phi_{0}$ is the base flag of $\cT$, label $\Phi_{0} \gamma $ with the triple $(\sigma, \vx, t)$. 
In particular, $\Phi_{0}$ is labelled with $((1), \cyvec{1}, 0)$ where $\cyvec{1}$ is the vector of $\cyctwo$ with all its entries equal to $1$.
Observe that we are just identifying every flag of $\cT$ with an element of its automorphism group. 
This identification is well defined, since $\aut(\cT)$ acts freely and transitively on flags. 
The combinatorial aspects of the triplets will be insightful in how the geometry and the combinatorics connect, as we shall explain below.

First observe that $\tras (\cT)$ acts on the facets of $\cT$ by translation. 
It follows that two flags  $(\sigma, \vx, t)$ and $(\tau, \vy, u)$ belong to the same facet of $\cT$ if an only if $u=t$. 
This allows us to identify the facets of $\cT$ with the elements of $\tras (\cT)$.

The first two coordinates of the label of a flag also have a combinatorial interpretation. 
We describe this interpretation only on the base facet $F_{0}$.
Recall that $F_{0}$ is an $n$-cube. 
Label the vertices of the cube with the elements of $\cyctwo$ in such a way that if $k \in \{0, \dots, n\}$, each $k$-face $F(\vx,I)$ of the cube can be described by its vertex set as \[F(\vx, I) =\{\vy  \in \cyctwo : y_{j} = x_{j} \text{ if } j \not\in I\},\] where $\vx = (x_{1}, \dots, x_{n})$, $\vy = (y_{1}, \dots, y_{n})$, $I \subset \{1, \dots, n\}$, $|I|=k$.  With this identification, two faces $F(\vx, I)$ and $G(\vx, J)$ with $|I| \leq |J|$ are incident if and only if $F(\vx,I) \subset G(\vy,J)$. In other words, we are just identifying the $n$-dimensional cube with the polytope $2^{\cS}$ where $\cS$ is the $(n-1)$-simplex (see \cite[Section 5.3]{Matousek_2002_LecturesDiscreteGeometry}).

Observe that  a flag of the cube containing a vertex $\vx$ now has the form  \[\{\emptyset, F(\vx, I_{0}), F(\vx, I_{1}), \dots , F(\vx, I_{n})\},\]  where $|I_{j}|= j$ for $j \in \{0, \dots, n\}$ and $I_{0} \subset I_{1} \subset \cdots \subset I_{n}$. Thus $I_{0}= \emptyset$ and $I_{n}= \left\{ 1 \dots, n \right\}$. Notice that for every $j \in \{1, \dots, n\}$, the set $I_{j} \sm I_{j-1}$ has exactly one element $i_{j}$.
The family of sets $\{I_{j} : 0 \leq j \leq n\}$ defines a permutation $\sigma \in S_{n}$ such that $ \sigma: j \mapsto i_{j}$. Conversely, a permutation $\sigma \in S_{n}$ determines a family $\{I_{j} : 0 \leq j \leq n\}$ of nested sets such that $I_{j} \sm I_{j-1} = j \sigma$.
Therefore, a permutation $\sigma$ and an element $\vx \in \cyctwo$ uniquely determine a flag of the base cube.
The label associated to this flag is precisely $( \sigma,\vx, 0)$.

\begin{figure}\centering
\def\svgwidth{.7\textwidth}
\begin{scriptsize}
	\begingroup \makeatletter \providecommand\color[2][]{\errmessage{(Inkscape) Color is used for the text in Inkscape, but the package 'color.sty' is not loaded}\renewcommand\color[2][]{}}\providecommand\transparent[1]{\errmessage{(Inkscape) Transparency is used (non-zero) for the text in Inkscape, but the package 'transparent.sty' is not loaded}\renewcommand\transparent[1]{}}\providecommand\rotatebox[2]{#2}\newcommand*\fsize{\dimexpr\f@size pt\relax}\newcommand*\lineheight[1]{\fontsize{\fsize}{#1\fsize}\selectfont}\ifx\svgwidth\undefined \setlength{\unitlength}{956.47699311bp}\ifx\svgscale\undefined \relax \else \setlength{\unitlength}{\unitlength * \real{\svgscale}}\fi \else \setlength{\unitlength}{\svgwidth}\fi \global\let\svgwidth\undefined \global\let\svgscale\undefined \makeatother \begin{picture}(1,0.98095816)\lineheight{1}\setlength\tabcolsep{0pt}\put(0,0){\includegraphics[width=\unitlength,page=1]{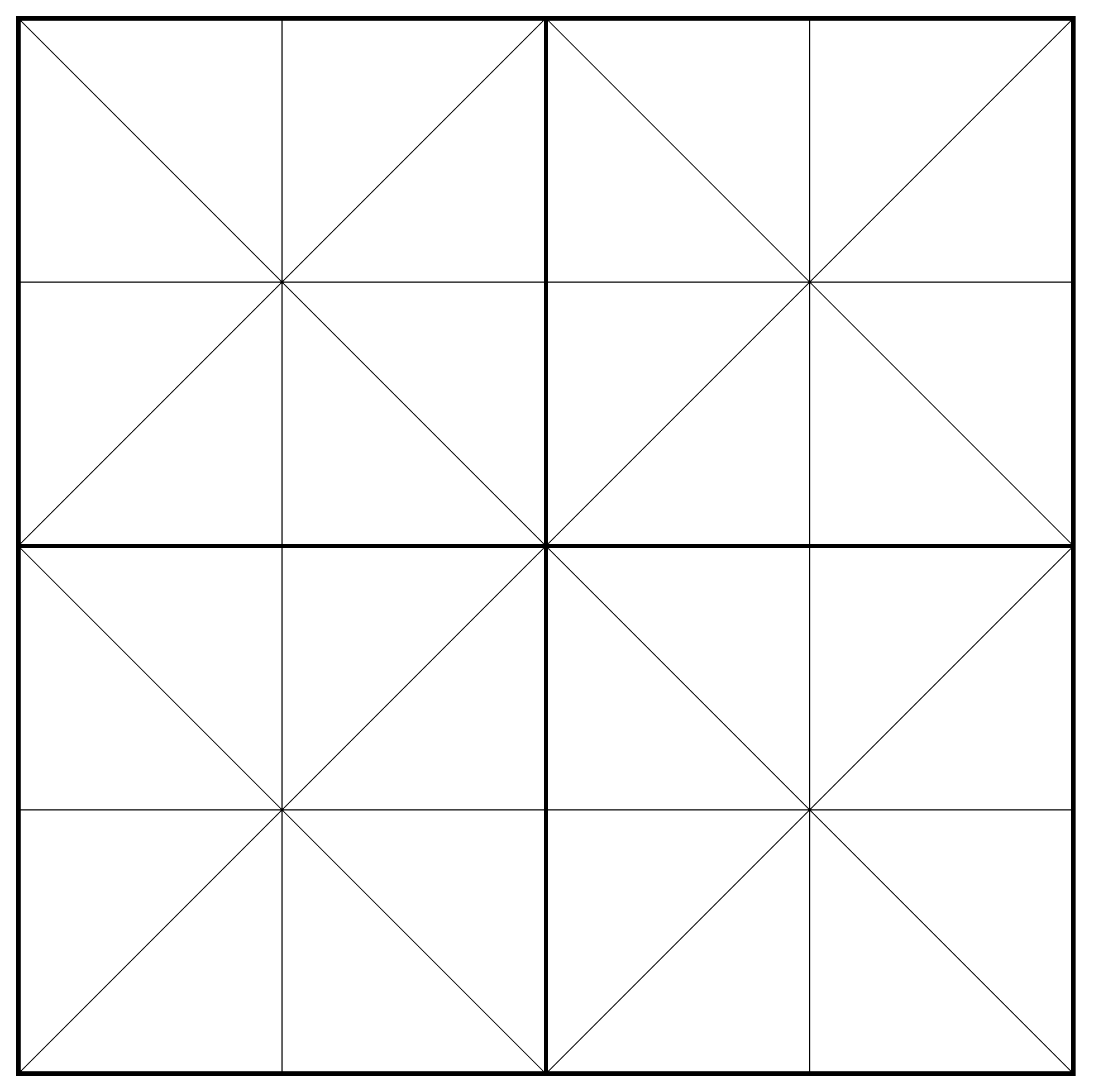}}\put(0.1848543,0.91512875){\color[rgb]{0,0,0}\makebox(0,0)[t]{\lineheight{1.25}\smash{\begin{tabular}[t]{c}$(1),$\end{tabular}}}}\put(0.1848543,0.88875711){\color[rgb]{0,0,0}\makebox(0,0)[t]{\lineheight{1.25}\smash{\begin{tabular}[t]{c}$(1,-1),$\end{tabular}}}}\put(0.18429418,0.86238546){\color[rgb]{0,0,0}\makebox(0,0)[t]{\lineheight{1.25}\smash{\begin{tabular}[t]{c}$0$\end{tabular}}}}\put(0.31821773,0.91512875){\color[rgb]{0,0,0}\makebox(0,0)[t]{\lineheight{1.25}\smash{\begin{tabular}[t]{c}$(1),$\end{tabular}}}}\put(0.31821773,0.88875711){\color[rgb]{0,0,0}\makebox(0,0)[t]{\lineheight{1.25}\smash{\begin{tabular}[t]{c}$(-1,-1),$\end{tabular}}}}\put(0.31765761,0.86238546){\color[rgb]{0,0,0}\makebox(0,0)[t]{\lineheight{1.25}\smash{\begin{tabular}[t]{c}$0$\end{tabular}}}}\put(0.41583966,0.81140164){\color[rgb]{0,0,0}\makebox(0,0)[t]{\lineheight{1.25}\smash{\begin{tabular}[t]{c}$(1\ 2),$\end{tabular}}}}\put(0.41583966,0.78502999){\color[rgb]{0,0,0}\makebox(0,0)[t]{\lineheight{1.25}\smash{\begin{tabular}[t]{c}$(-1,-1),$\end{tabular}}}}\put(0.41527954,0.75865834){\color[rgb]{0,0,0}\makebox(0,0)[t]{\lineheight{1.25}\smash{\begin{tabular}[t]{c}$0$\end{tabular}}}}\put(0.41583966,0.6780382){\color[rgb]{0,0,0}\makebox(0,0)[t]{\lineheight{1.25}\smash{\begin{tabular}[t]{c}$(1\ 2),$\end{tabular}}}}\put(0.41583966,0.65166656){\color[rgb]{0,0,0}\makebox(0,0)[t]{\lineheight{1.25}\smash{\begin{tabular}[t]{c}$(-1,1),$\end{tabular}}}}\put(0.41527954,0.62529491){\color[rgb]{0,0,0}\makebox(0,0)[t]{\lineheight{1.25}\smash{\begin{tabular}[t]{c}$0$\end{tabular}}}}\put(0.33303589,0.60394741){\color[rgb]{0,0,0}\makebox(0,0)[t]{\lineheight{1.25}\smash{\begin{tabular}[t]{c}$(1),$\end{tabular}}}}\put(0.33303589,0.57757577){\color[rgb]{0,0,0}\makebox(0,0)[t]{\lineheight{1.25}\smash{\begin{tabular}[t]{c}$(-1,1),$\end{tabular}}}}\put(0.33247577,0.55120412){\color[rgb]{0,0,0}\makebox(0,0)[t]{\lineheight{1.25}\smash{\begin{tabular}[t]{c}$0$\end{tabular}}}}\put(0.18485429,0.60394741){\color[rgb]{0,0,0}\makebox(0,0)[t]{\lineheight{1.25}\smash{\begin{tabular}[t]{c}$(1),$\end{tabular}}}}\put(0.18485429,0.57757577){\color[rgb]{0,0,0}\makebox(0,0)[t]{\lineheight{1.25}\smash{\begin{tabular}[t]{c}$(1,1),$\end{tabular}}}}\put(0.18429418,0.55120412){\color[rgb]{0,0,0}\makebox(0,0)[t]{\lineheight{1.25}\smash{\begin{tabular}[t]{c}$0$\end{tabular}}}}\put(0.08112718,0.6780382){\color[rgb]{0,0,0}\makebox(0,0)[t]{\lineheight{1.25}\smash{\begin{tabular}[t]{c}$(1\ 2),$\end{tabular}}}}\put(0.08112718,0.65166656){\color[rgb]{0,0,0}\makebox(0,0)[t]{\lineheight{1.25}\smash{\begin{tabular}[t]{c}$(1,1),$\end{tabular}}}}\put(0.08056707,0.62529491){\color[rgb]{0,0,0}\makebox(0,0)[t]{\lineheight{1.25}\smash{\begin{tabular}[t]{c}$0$\end{tabular}}}}\put(0.08112718,0.81140164){\color[rgb]{0,0,0}\makebox(0,0)[t]{\lineheight{1.25}\smash{\begin{tabular}[t]{c}$(1\ 2),$\end{tabular}}}}\put(0.08112718,0.78502999){\color[rgb]{0,0,0}\makebox(0,0)[t]{\lineheight{1.25}\smash{\begin{tabular}[t]{c}$(1,-1),$\end{tabular}}}}\put(0.08056707,0.75865834){\color[rgb]{0,0,0}\makebox(0,0)[t]{\lineheight{1.25}\smash{\begin{tabular}[t]{c}$0$\end{tabular}}}}\put(0.665652,0.91512875){\color[rgb]{0,0,0}\makebox(0,0)[t]{\lineheight{1.25}\smash{\begin{tabular}[t]{c}$(1),$\end{tabular}}}}\put(0.665652,0.88875711){\color[rgb]{0,0,0}\makebox(0,0)[t]{\lineheight{1.25}\smash{\begin{tabular}[t]{c}$(1,-1),$\end{tabular}}}}\put(0.66509188,0.86238546){\color[rgb]{0,0,0}\makebox(0,0)[t]{\lineheight{1.25}\smash{\begin{tabular}[t]{c}$e_1$\end{tabular}}}}\put(0.79901543,0.91512875){\color[rgb]{0,0,0}\makebox(0,0)[t]{\lineheight{1.25}\smash{\begin{tabular}[t]{c}$(1),$\end{tabular}}}}\put(0.79901543,0.88875711){\color[rgb]{0,0,0}\makebox(0,0)[t]{\lineheight{1.25}\smash{\begin{tabular}[t]{c}$(-1,-1),$\end{tabular}}}}\put(0.79845532,0.86238546){\color[rgb]{0,0,0}\makebox(0,0)[t]{\lineheight{1.25}\smash{\begin{tabular}[t]{c}$e_1$\end{tabular}}}}\put(0.89663737,0.81140164){\color[rgb]{0,0,0}\makebox(0,0)[t]{\lineheight{1.25}\smash{\begin{tabular}[t]{c}$(1\ 2),$\end{tabular}}}}\put(0.89663737,0.78502999){\color[rgb]{0,0,0}\makebox(0,0)[t]{\lineheight{1.25}\smash{\begin{tabular}[t]{c}$(-1,-1),$\end{tabular}}}}\put(0.89607726,0.75865834){\color[rgb]{0,0,0}\makebox(0,0)[t]{\lineheight{1.25}\smash{\begin{tabular}[t]{c}$e_1$\end{tabular}}}}\put(0.89663737,0.6780382){\color[rgb]{0,0,0}\makebox(0,0)[t]{\lineheight{1.25}\smash{\begin{tabular}[t]{c}$(1\ 2),$\end{tabular}}}}\put(0.89663737,0.65166656){\color[rgb]{0,0,0}\makebox(0,0)[t]{\lineheight{1.25}\smash{\begin{tabular}[t]{c}$(-1,1),$\end{tabular}}}}\put(0.89607726,0.62529491){\color[rgb]{0,0,0}\makebox(0,0)[t]{\lineheight{1.25}\smash{\begin{tabular}[t]{c}$e_1$\end{tabular}}}}\put(0.81383359,0.60394741){\color[rgb]{0,0,0}\makebox(0,0)[t]{\lineheight{1.25}\smash{\begin{tabular}[t]{c}$(1),$\end{tabular}}}}\put(0.81383359,0.57757577){\color[rgb]{0,0,0}\makebox(0,0)[t]{\lineheight{1.25}\smash{\begin{tabular}[t]{c}$(-1,1),$\end{tabular}}}}\put(0.81327347,0.55120412){\color[rgb]{0,0,0}\makebox(0,0)[t]{\lineheight{1.25}\smash{\begin{tabular}[t]{c}$e_1$\end{tabular}}}}\put(0.665652,0.60394741){\color[rgb]{0,0,0}\makebox(0,0)[t]{\lineheight{1.25}\smash{\begin{tabular}[t]{c}$(1),$\end{tabular}}}}\put(0.665652,0.57757577){\color[rgb]{0,0,0}\makebox(0,0)[t]{\lineheight{1.25}\smash{\begin{tabular}[t]{c}$(1,1),$\end{tabular}}}}\put(0.66509188,0.55120412){\color[rgb]{0,0,0}\makebox(0,0)[t]{\lineheight{1.25}\smash{\begin{tabular}[t]{c}$e_1$\end{tabular}}}}\put(0.56192488,0.6780382){\color[rgb]{0,0,0}\makebox(0,0)[t]{\lineheight{1.25}\smash{\begin{tabular}[t]{c}$(1\ 2),$\end{tabular}}}}\put(0.56192488,0.65166656){\color[rgb]{0,0,0}\makebox(0,0)[t]{\lineheight{1.25}\smash{\begin{tabular}[t]{c}$(1,1),$\end{tabular}}}}\put(0.56136477,0.62529491){\color[rgb]{0,0,0}\makebox(0,0)[t]{\lineheight{1.25}\smash{\begin{tabular}[t]{c}$e_1$\end{tabular}}}}\put(0.56192488,0.81140164){\color[rgb]{0,0,0}\makebox(0,0)[t]{\lineheight{1.25}\smash{\begin{tabular}[t]{c}$(1\ 2),$\end{tabular}}}}\put(0.56192488,0.78502999){\color[rgb]{0,0,0}\makebox(0,0)[t]{\lineheight{1.25}\smash{\begin{tabular}[t]{c}$(1,-1),$\end{tabular}}}}\put(0.56136477,0.75865834){\color[rgb]{0,0,0}\makebox(0,0)[t]{\lineheight{1.25}\smash{\begin{tabular}[t]{c}$e_1$\end{tabular}}}}\put(0.1848543,0.4367073){\color[rgb]{0,0,0}\makebox(0,0)[t]{\lineheight{1.25}\smash{\begin{tabular}[t]{c}$(1),$\end{tabular}}}}\put(0.1848543,0.41033565){\color[rgb]{0,0,0}\makebox(0,0)[t]{\lineheight{1.25}\smash{\begin{tabular}[t]{c}$(1,-1),$\end{tabular}}}}\put(0.18429418,0.383964){\color[rgb]{0,0,0}\makebox(0,0)[t]{\lineheight{1.25}\smash{\begin{tabular}[t]{c}$e_2$\end{tabular}}}}\put(0.31821773,0.4367073){\color[rgb]{0,0,0}\makebox(0,0)[t]{\lineheight{1.25}\smash{\begin{tabular}[t]{c}$(1),$\end{tabular}}}}\put(0.31821773,0.41033565){\color[rgb]{0,0,0}\makebox(0,0)[t]{\lineheight{1.25}\smash{\begin{tabular}[t]{c}$(-1,-1),$\end{tabular}}}}\put(0.31765761,0.383964){\color[rgb]{0,0,0}\makebox(0,0)[t]{\lineheight{1.25}\smash{\begin{tabular}[t]{c}$e_2$\end{tabular}}}}\put(0.41583966,0.33298019){\color[rgb]{0,0,0}\makebox(0,0)[t]{\lineheight{1.25}\smash{\begin{tabular}[t]{c}$(1\ 2),$\end{tabular}}}}\put(0.41583966,0.30660854){\color[rgb]{0,0,0}\makebox(0,0)[t]{\lineheight{1.25}\smash{\begin{tabular}[t]{c}$(-1,-1),$\end{tabular}}}}\put(0.41527954,0.28023689){\color[rgb]{0,0,0}\makebox(0,0)[t]{\lineheight{1.25}\smash{\begin{tabular}[t]{c}$e_2$\end{tabular}}}}\put(0.41583966,0.19961675){\color[rgb]{0,0,0}\makebox(0,0)[t]{\lineheight{1.25}\smash{\begin{tabular}[t]{c}$(1\ 2),$\end{tabular}}}}\put(0.41583966,0.1732451){\color[rgb]{0,0,0}\makebox(0,0)[t]{\lineheight{1.25}\smash{\begin{tabular}[t]{c}$(-1,1),$\end{tabular}}}}\put(0.41527954,0.14687346){\color[rgb]{0,0,0}\makebox(0,0)[t]{\lineheight{1.25}\smash{\begin{tabular}[t]{c}$e_2$\end{tabular}}}}\put(0.33303589,0.12552599){\color[rgb]{0,0,0}\makebox(0,0)[t]{\lineheight{1.25}\smash{\begin{tabular}[t]{c}$(1),$\end{tabular}}}}\put(0.33303589,0.09915434){\color[rgb]{0,0,0}\makebox(0,0)[t]{\lineheight{1.25}\smash{\begin{tabular}[t]{c}$(-1,1),$\end{tabular}}}}\put(0.33247577,0.07278269){\color[rgb]{0,0,0}\makebox(0,0)[t]{\lineheight{1.25}\smash{\begin{tabular}[t]{c}$e_2$\end{tabular}}}}\put(0.18485429,0.12552599){\color[rgb]{0,0,0}\makebox(0,0)[t]{\lineheight{1.25}\smash{\begin{tabular}[t]{c}$(1),$\end{tabular}}}}\put(0.18485429,0.09915434){\color[rgb]{0,0,0}\makebox(0,0)[t]{\lineheight{1.25}\smash{\begin{tabular}[t]{c}$(1,1),$\end{tabular}}}}\put(0.18429418,0.07278269){\color[rgb]{0,0,0}\makebox(0,0)[t]{\lineheight{1.25}\smash{\begin{tabular}[t]{c}$e_2$\end{tabular}}}}\put(0.08112718,0.19961675){\color[rgb]{0,0,0}\makebox(0,0)[t]{\lineheight{1.25}\smash{\begin{tabular}[t]{c}$(1\ 2),$\end{tabular}}}}\put(0.08112718,0.1732451){\color[rgb]{0,0,0}\makebox(0,0)[t]{\lineheight{1.25}\smash{\begin{tabular}[t]{c}$(1,1),$\end{tabular}}}}\put(0.08056707,0.14687346){\color[rgb]{0,0,0}\makebox(0,0)[t]{\lineheight{1.25}\smash{\begin{tabular}[t]{c}$e_2$\end{tabular}}}}\put(0.08112718,0.33298019){\color[rgb]{0,0,0}\makebox(0,0)[t]{\lineheight{1.25}\smash{\begin{tabular}[t]{c}$(1\ 2),$\end{tabular}}}}\put(0.08112718,0.30660854){\color[rgb]{0,0,0}\makebox(0,0)[t]{\lineheight{1.25}\smash{\begin{tabular}[t]{c}$(1,-1),$\end{tabular}}}}\put(0.08056707,0.28023689){\color[rgb]{0,0,0}\makebox(0,0)[t]{\lineheight{1.25}\smash{\begin{tabular}[t]{c}$e_2$\end{tabular}}}}\put(0.66644409,0.4367073){\color[rgb]{0,0,0}\makebox(0,0)[t]{\lineheight{1.25}\smash{\begin{tabular}[t]{c}$(1),$\end{tabular}}}}\put(0.66644409,0.41033565){\color[rgb]{0,0,0}\makebox(0,0)[t]{\lineheight{1.25}\smash{\begin{tabular}[t]{c}$(1,-1),$\end{tabular}}}}\put(0.66588397,0.383964){\color[rgb]{0,0,0}\makebox(0,0)[t]{\lineheight{1.25}\smash{\begin{tabular}[t]{c}$e_1 +e_2$\end{tabular}}}}\put(0.79980752,0.4367073){\color[rgb]{0,0,0}\makebox(0,0)[t]{\lineheight{1.25}\smash{\begin{tabular}[t]{c}$(1),$\end{tabular}}}}\put(0.79980752,0.41033565){\color[rgb]{0,0,0}\makebox(0,0)[t]{\lineheight{1.25}\smash{\begin{tabular}[t]{c}$(-1,-1),$\end{tabular}}}}\put(0.79924741,0.383964){\color[rgb]{0,0,0}\makebox(0,0)[t]{\lineheight{1.25}\smash{\begin{tabular}[t]{c}$e_1 +e_2$\end{tabular}}}}\put(0.89742946,0.33298019){\color[rgb]{0,0,0}\makebox(0,0)[t]{\lineheight{1.25}\smash{\begin{tabular}[t]{c}$(1\ 2),$\end{tabular}}}}\put(0.89742946,0.30660854){\color[rgb]{0,0,0}\makebox(0,0)[t]{\lineheight{1.25}\smash{\begin{tabular}[t]{c}$(-1,-1),$\end{tabular}}}}\put(0.89686935,0.28023689){\color[rgb]{0,0,0}\makebox(0,0)[t]{\lineheight{1.25}\smash{\begin{tabular}[t]{c}$e_1 +e_2$\end{tabular}}}}\put(0.89742946,0.19961675){\color[rgb]{0,0,0}\makebox(0,0)[t]{\lineheight{1.25}\smash{\begin{tabular}[t]{c}$(1\ 2),$\end{tabular}}}}\put(0.89742946,0.1732451){\color[rgb]{0,0,0}\makebox(0,0)[t]{\lineheight{1.25}\smash{\begin{tabular}[t]{c}$(-1,1),$\end{tabular}}}}\put(0.89686935,0.14687346){\color[rgb]{0,0,0}\makebox(0,0)[t]{\lineheight{1.25}\smash{\begin{tabular}[t]{c}$e_1 +e_2$\end{tabular}}}}\put(0.81462568,0.12552599){\color[rgb]{0,0,0}\makebox(0,0)[t]{\lineheight{1.25}\smash{\begin{tabular}[t]{c}$(1),$\end{tabular}}}}\put(0.81462568,0.09915434){\color[rgb]{0,0,0}\makebox(0,0)[t]{\lineheight{1.25}\smash{\begin{tabular}[t]{c}$(-1,1),$\end{tabular}}}}\put(0.81406556,0.07278269){\color[rgb]{0,0,0}\makebox(0,0)[t]{\lineheight{1.25}\smash{\begin{tabular}[t]{c}$e_1 +e_2$\end{tabular}}}}\put(0.66644409,0.12552599){\color[rgb]{0,0,0}\makebox(0,0)[t]{\lineheight{1.25}\smash{\begin{tabular}[t]{c}$(1),$\end{tabular}}}}\put(0.66644409,0.09915434){\color[rgb]{0,0,0}\makebox(0,0)[t]{\lineheight{1.25}\smash{\begin{tabular}[t]{c}$(1,1),$\end{tabular}}}}\put(0.66588397,0.07278269){\color[rgb]{0,0,0}\makebox(0,0)[t]{\lineheight{1.25}\smash{\begin{tabular}[t]{c}$e_1 +e_2$\end{tabular}}}}\put(0.56271697,0.19961675){\color[rgb]{0,0,0}\makebox(0,0)[t]{\lineheight{1.25}\smash{\begin{tabular}[t]{c}$(1\ 2),$\end{tabular}}}}\put(0.56271697,0.1732451){\color[rgb]{0,0,0}\makebox(0,0)[t]{\lineheight{1.25}\smash{\begin{tabular}[t]{c}$(1,1),$\end{tabular}}}}\put(0.56215686,0.14687346){\color[rgb]{0,0,0}\makebox(0,0)[t]{\lineheight{1.25}\smash{\begin{tabular}[t]{c}$e_1 +e_2$\end{tabular}}}}\put(0.56271697,0.33298019){\color[rgb]{0,0,0}\makebox(0,0)[t]{\lineheight{1.25}\smash{\begin{tabular}[t]{c}$(1\ 2),$\end{tabular}}}}\put(0.56271697,0.30660854){\color[rgb]{0,0,0}\makebox(0,0)[t]{\lineheight{1.25}\smash{\begin{tabular}[t]{c}$(1,-1),$\end{tabular}}}}\put(0.56215686,0.28023689){\color[rgb]{0,0,0}\makebox(0,0)[t]{\lineheight{1.25}\smash{\begin{tabular}[t]{c}$e_1 +e_2$\end{tabular}}}}\put(0,0){\includegraphics[width=\unitlength,page=2]{chirExtToroids_Labeling.pdf}}\end{picture}\endgroup  \end{scriptsize}
\caption{Labeling of the flags of the toroid $\{4,4\}_{(2,0)}$}\label{fig:labeling}
\end{figure}

Informally speaking, if $(\sigma, \vx, 0)$ is the label of a flag $\Phi$, then $\vx$ describes the relative position of the vertex of $\Phi$ on the cube $F_{0}$. The permutation $\sigma$ defines the ``direction'' of the faces relative to $\vx$ in the following sense. To get the other vertex of the edge in $\Phi$ we have to ``move'' (change the sign of $\vx$) in direction $1 \sigma$; to get the four vertices of the $2$-face of $\Phi$ we have to ``allow movement'' in the directions $1 \sigma$ and $2 \sigma$; etc. See Figure\nobreakspace \ref {fig:labeling} for an example on dimension two.

We next introduce some notation. Given $\vx = (x_{1}, \dots, x_{n}) \in \cyctwo$ and a subset $I=\{i_{1},\dots, i_{k}\} $ of $ \{1, \dots, n\}$ we denote by $\vx^{(i_{1}, \dots, i_{k})}$ the vector $(x'_{1}, \dots, x'_{n}) \in \cyctwo$ such that $x_{i} = x'_{i}$ if and only if $i \not\in \{i_{1}, \dots, i_{k}\}$. Note that $\vx^{(i_{1}, \dots, i_{k})} = \vy_{I} \vx$ where $\vy_{I}$ is the vector in $\cyctwo$ whose $i^{th}$-entry is $1$ if and only if $i \not \in I$. 
For a permutation $\sigma \in S_{n}$ and a vector $\vx = (x_{1}, \dots, x_{n}) \in \cyctwo$, we denote by $\sigma \vx$ the vector resulting after permuting the coordinates of $\vx$ according to $\sigma$, this is, $\sigma \vx = (x_{1 \sigma^{-1}}, x_{2 \sigma^{-1}}, \dots, x_{n \sigma^{-1}} )$. 
Similarly, if $t= (t_{1}, \dots, t_{n}) \in \tras (\cT)$, then $\sigma t$ denotes the vector $(t_{1 \sigma^{-1}}, t_{2 \sigma^{-1}}, \dots, t_{n \sigma^{-1}} )$, and if $\vy =(y_{1}, \dots, y_{n}) \in \cyctwo$, then $\vy t =(y_{1}t_{1}, \dots, y_{n} t_{n})$. 
Finally, observe that all these operations can be understood as left actions of the corresponding groups.

With the notation just introduced it is easy to describe the action of $\mon(\cT)$ on the set of flags. This is given by

\begin{equation} \label{eq:monOnLabels}
	r_{i}(\sigma, \vx,  t) =
	\begin{cases}
		 ( \sigma, \vx^{(1\sigma)}, t) & \text{if}\  i=0, \\
		 ((i\ i+1) \sigma, \vx,  t) & \text{if}\ 1 \leq i \leq n-1, \\
		 ( \sigma, \vx^{(n \sigma)}, t-x_{n\sigma}e_{n\sigma}) & \text{if}\ i=n,\\
	\end{cases}
\end{equation}
where $\vx =  (x_{1}, \dots, x_{n})$ and $e_{k}$ denotes the vector $(0^{k-1}, 1, 0^{n-k}) \in\tras (\cT)$. 

The validity of Equation\nobreakspace \textup {(\ref {eq:monOnLabels})} can be derived from the combinatorial interpretation of the labelling discussed above. 
For instance, it is clear that $r_{0}$ must fix the first and third coordinates in $\Phi = (\sigma, \vx, t)$ since the faces of rank $i\geq 1$ of $\Phi$ and $r_0 \Phi$ are the same, but the vertex of $r_0 \Phi$ is a vertex that result from moving $\vx$ in the direction of the base edge ($1 \sigma$).
Alternatively, we can also define the $(n+1)$-maniplex $\cM$ (in the sense of \cite{Wilson_2012_ManiplexesPart1}) whose flags are the set $S_n \times C_{2}^{n} \times \tras(\cT)$ and whose monodromy group is given by Equation\nobreakspace \textup {(\ref {eq:monOnLabels})} and prove that $\cM \cong \cT$.
We leave the details to the reader. 

Assume now that $ \gamma \in \aut(\cT) = (S_{n} \ltimes \cyctwo) \ltimes \tras (\cT)$ and let $\gamma = \tau \cdot \vy \cdot u$, $\tau \in S_n$, $\vy \in \cyctwo$, $u \in \tras(\cT)$.
Let $\Phi = (\sigma, \vx, t)$ be a flag. 
Observe that $\tau^{-1} \vx \tau \in \cyctwo$ and that $\tau^{-1} \vx \tau$ is the automorphism given by the vector $\tau \vx$. 
Similarly, $\tau^{-1} t \tau$ is the automorphism given by the vector $\tau t \in \tras (\cT)$. 
For $\vy \in \cyctwo$, the automorphism $\vy^{-1} t \vy$ is given by the vector $\vy t \in \tras (\cT) $. 
The previous discussion implies that the action of $\aut(\cT)$ is given by
\begin{equation}\label{eq:autOnLabels_s}
  \begin{aligned}
    (\sigma, \vx, t) \tau &= (\sigma \tau, \tau \vx, \tau t),\\
    (\sigma, \vx, t) \vy &= (\sigma, \vy \vx, \vy t), \\
    (\sigma, \vx, t) u &= (\sigma, \vx, t + u), \\
  \end{aligned}
\end{equation}
or equivalently 
\begin{equation}\label{eq:autOnLabels}
(\sigma, \vx, t)\gamma = (\sigma \tau , \vy(\tau \vx),\vy(\tau t) + u ).
\end{equation}

As in Section\nobreakspace \ref {sec:rotary}, we  say that a flag $\Psi$ of $\cT$ is \emph{white} if $\Psi$  and the base flag $\Phi_{0}$ belong to the same orbit under $\autp(\cT)$. 
If $\Phi$ is not white, then we say that $\Phi$ is a \emph{black} flag.  
If $\sigma \in S_{n}$, then $\sgn(\sigma)\in \{1,-1\}$ and is equal to $1$ if and only if $\sigma$ is an even permutation. 
If $\vx \in \cyctwo$, let $\sgn(\vx) = (-1)^{k}$ where $k$ denotes the number of entries of $\vx$ equal to $-1$. 
Observe that a flag $(\sigma, \vx, t)$ of $\cT$ is white if and only if $\sgn(\sigma)\sgn(\vx) = 1$.

Let us now prove some structural results related to the labelling presented above that will be of use in Section\nobreakspace \ref {sec:chirExtRegToroids}.

\begin{lem}\label{lem:flagsAroundX0} 
  Let $X_{0}$ be the base vertex of $\cT$. 
  Let $\Phi = (\sigma, \vx, t) \in  \cF(\cT)$ with $\vx=(x_{1}, \dots, x_{n})$ and $t= (t_{1}, \dots, t_{n})$. 
  Then $\Phi$ contains $X_{0}$ if and only if $t_i \in \{0,-1\} $ and $t_i = -1$ if and only if $x_i = -1$.
\end{lem}
  \begin{proof}
    Assume first that $\Phi$ is a flag containing the vertex $X_0$.
    Let $\Phi_0 = ( (1), \cyvec{1},0)$ be the base flag of $\cT$. 
    Since  $\Phi$ and $\Phi_0$ share the vertex, then there exists $w \in \langle r_1, \dots, r_n \rangle$ such that $\Phi= w \Phi_0$.
    Let $v_0, \dots, v_d \in \langle r_1, \dots, r_{n-1} \rangle$ be such that \[w= v_d r_n v_{d-1} \cdots v_1 r_n v_0. \]
    Observe that $v_i$ (for $ 1 \leq i \leq d )$ preserves the second and third entries of the label of every flag. Meanwhile $r_n$ changes both the second and third entries, and it changes them on the same coordinate.
    Since both such entries of the label of $\Phi_0$ are trivial, an inductive argument on $d$ shows that $t_i \in \{0, -1\}$ and $t_i = -1$ if and only if $x_i = -1$.
  \end{proof}

An immediate consequence of Lemma\nobreakspace \ref {lem:flagsAroundX0} is the following result.

\begin{coro}
\label{coro:facetsAroundAnyVertex}
Let $\Phi$ be the flag $(\sigma, \vx, t)$, with $\vx=(x_{1}, \dots, x_{n})$ and $t= (t_{1}, \dots, t_{n})$. 
Let $X$ be a vertex of $\cT$ and $u = (u_1, \dots, u_n) \in \tras(\cT)$ be the (unique) translation of $\cT$ that maps the base vertex $X_0$ to $X$.
Then $\Phi$ contains the vertex $X$ if and only if $t_i-u_i \in \{0, -1\}$ and $t_i -u_i = -1$ if and only if $x_i=-1$.
In particular, the facets of $\cT$ that contain the vertex $X$ are those whose associated vector is in the set \[\left\{(u_1 + \epsilon_1, \dots, u_n + \epsilon_n) : \epsilon_1, \dots \epsilon_n \in \{0,-1\} \right\}.\]
\end{coro}

 \section{Chiral extensions of regular toroids}\label{sec:chirExtRegToroids}

In this section we develop the construction of (the automorphism group of) a chiral extension of a regular toroid $\cT$ satisfying certain properties (described below).
As a consequence we prove Theorem\nobreakspace \ref {thm:elChido}.

The construction of a chiral extension of a toroid $\cT$ relies on certain properties of a particular monodromy element $\kng$.
We describe such properties in the following paragraphs.

Let $((1), \cyvec{1},0)$ be the base flag of $\cT$.
Let $v_{0}$ 
be the vector $(1,2,\dots,n) \in \tras(\cT)$.
Let $\kng$ be the (unique) element of $\monp(\cT)$ that satisfies that $\kng ((1), \cyvec{1},0) = ( (1), \cyvec{1}, v_{0})$.

For a permutation $\sigma \in S_{n}$ and an element $\vx \in \cyctwo$, let $v(\sigma, \vx)$ denote the vector $\vx \sigma  v_{0}$. This is the vector obtained from $v_{0}$ after permuting its entries according to $\sigma$ then changing signs according to $\vx$. It follows from Equation\nobreakspace \textup {(\ref {eq:autOnLabels})} that if $(\sigma, \vx,  t)$ is a flag, then
\begin{equation} \label{eq:eta}
  \kng(\sigma, \vx,  t) = \left(\sigma, \vx,  t + v(\sigma, \vx)\right).
\end{equation}
For a flag $\Phi = (\sigma, \vx, t)$ of $\cT$, let $\hat{\Phi}$ denote the flag $(\sigma, -\vx, t)$.
In this particular case Equation\nobreakspace \textup {(\ref {eq:eta})} takes the form of
\begin{equation}\label{eq:antipoda}
\kng \hat{\Phi} = (\sigma, \vx, t - v(\sigma, \vx)).
\end{equation}

Consider also the element $\bar{\kng} = r_{0} \kng r_{0}$.
Observe that
\begin{equation} \label{eq:etaBar}
\bar{\kng}(\sigma, \vx,  t) = \left(\sigma, \vx,  t + v(\sigma, \vx^{1 \sigma})\right).
\end{equation}

\begin{rem}\label{rem:largeToroid} 
In what follows we will use the permutations $\kng$ and $\bar{\kng}$ to define a set of special facets of $\cT$. 
We will do so by considering the images of the base facet of $\cT$ under some power of $\kng$ or $\bar{\kng}$. 
Recall that $\cT$ is defined by a vector of the form $(a,0,\ldots,0)$, $(a,a,0,\ldots,0)$ or $(a,a,\ldots,a)$ for some integer $a$. 
For things to go smoothly, and to ensure that all such special facets are distinct from one another, we need $\cT$ to be `large enough'.
More precisely, we shall require that $a \geq 6n + 1$.
This is a technical requirement we impose for the sake of convenience, to ensure that $\cT$ is large enough so that we do not end up `looping around' the toroid when we do not want to.
Moreover, when talking about translations in $\tras(\cT)$ we can think of their translation vectors as elements in $\bZ^n$ and not in $\bZ^n/\LL$, or more precisely, as their representatives in the appropriate $D(\bLL)$.
We will henceforth always assume that the parameter $a$ satisfies the above bounds, even though this might not be strictly necessary for the construction we are about to describe to work.
\end{rem}

 For a flag $\Psi \in \cF(\cT)$, let $F(\Psi)$ denote the facet of $\Psi$.
 If $\Psi = (\sigma, \vx, t)$ we shall abuse notation and use $F(\Psi)$ to denote both the facet in the poset of $\cT$ and the vector $t$.

\begin{prop} \label{prop:etaIsGood}
Let $F_{0}$ be the base facet of $\cT$.
For $i,j \in \{\pm1, \pm 2, \pm 3\}$ and $\Psi_{1} \neq \Psi_{2}$ flags of $\cT$ containing the facet $F$ of $\cT$. 
The elements $\kng$ and $\etab$ satisfy:
\begin{enumerate}
    \item \label{item:baseFacet} If $F = F_{0}$, then the facets $F(\kng^{i} \Psi_{1})$ and $F(\kng^{i} \Psi_{2})$ are different.
    Similarly, the facets  $F(\etab^{i} \Psi_{1})$ and $F(\etab^{i} \Psi_{2})$ are different.
    \item \label{item:anyFacet} If $F$ is any arbitrary facet, then $F(\kng^{i} \Psi_{1}) \neq F(\kng^{i} \Psi_{2})$.
    Similarly $F(\etab^{i} \Psi_{1}) \neq F(\etab^{i} \Psi_{2})$.
    \item \label{item:barEta} The flags $\kng \Psi_{1}$ and $\etab \Psi_{1}^{0}=r_0 \kng r_0 r_0 \Psi_{1}$ are $0$-adjacent.
    In particular, $F(\kng \Psi_{1}) = F(\etab\Psi_{1}^{0})$.
\item \label{item:Antipoda} The facet of $\kng^{i} \Psi_{1}$ is the same as the facet of $\kng^{-i} \hat{\Psi_{1}} $.
    Similarly $F(\etab^{i} \Psi_{1}) = F(\etab^{-i} \hat{\Psi_{1}})$.
    \item \label{item:NMA} If $i \neq \pm j$, then $\left\{F(\kng^{i} \Psi_{1}), F(\etab^{i} \Psi_{1})\right\} \cap \left\{ F(\kng^{j} \Psi_{2}), F(\etab^{j}\Psi_{2}) \right\} = \emptyset$.
    \item \label{item:esAntipoda} If 
either $F(\kng^{i} \Psi_{1}) = F(\kng^{j} \Psi_{2}) $ or $F(\etab^{i} \Psi_{1}) = F(\etab^{j} \Psi_{2}) $, then $i = -j$ and $\Psi_{2}=\hat{\Psi_{1}}$.
    \item \label{item:esVecina} If 
$F(\kng^{i} \Psi_{1}) = F(\etab^{j} \Psi_{2}) $, then $i=j$ and $\Psi_{1} = r_{0}\Psi_{2}$; or $i = -j$ and $\Psi_{1}=r_0 \hat{\Psi_{2}}$.

\end{enumerate}
\end{prop}
\begin{proof}
    If $\Psi_1$ and $\Psi_2$ are different flags containing the base facet $F_0$, then $\Psi_1 = (\sigma_{1}, \vx_{1}, 0)$ and $\Psi_{2} = (\sigma_{2}, \vx_{2}, 0)$ for some $\sigma_1, \sigma_{2} \in S_{n}$ and $\vx_{1}, \vx_{2} \in \cyctwo$ with $(\sigma_{1}, \vx_{1}) \neq (\sigma_{2}, \vx_{2})$.
    It follows that $\kng^{i}\Psi_1 = (\sigma_{1}, \vx_{1}, i v(\sigma_{1},\vx_1))$ and
    $\kng^{i}\Psi_{2} = (\sigma_{2}, \vx_{2}, i v(\sigma_{2},\vx_2))$.
    From our assumptions on $\cT$, the elements  $i v(\sigma_{1},\vx_1), i v(\sigma_{2},\vx_2) \in \tras(\cT)$ are different whenever $(\sigma_{1}, \vx_{1}) \neq (\sigma_{2}, \vx_{2})$, therefore, the facets of $\kng^{i}\Psi_{1}$ and $\kng^{i}\Psi_{2}$ are different.
    An analogous argument for $\etab$ completes the proof of Item\nobreakspace \ref {item:baseFacet}.
    Item\nobreakspace \ref {item:anyFacet} follows from Item\nobreakspace \ref {item:baseFacet} and from the fact that $\cT$ is facet-transitive. 
    In fact, assume that $G_1$ and $G_2$ are the facets of $\kng^{i} \Psi_1$ and $\kng^{i}\Psi_2$, respectively.
    Take $\gamma \in \aut(\cT)$ such that $F\gamma = F_0$. 
    The facet of $(\kng^{i} \Psi_1)\gamma = \kng^{i} (\Psi_1 \gamma)$ is $G_1 \gamma$, while $G_2 \gamma$ is the facet of $(\kng^{i} \Psi_2)\gamma = \kng^{i} (\Psi_2 \gamma)$. 
    Both $\Psi_1 \gamma$ and $\Psi_2 \gamma$ are flags containing $F_0$, hence $G_1 \gamma \neq G_2 \gamma$, which implies that $G_1 \neq G_2$.

    Item\nobreakspace \ref {item:barEta} is obvious.  Item\nobreakspace \ref {item:Antipoda} follows from Equation\nobreakspace \textup {(\ref {eq:antipoda})}.

    Let us prove Item\nobreakspace \ref {item:NMA} assuming that $F = F_{0}$.
    The other cases follow in a similar fashion as Item\nobreakspace \ref {item:anyFacet} follows from Item\nobreakspace \ref {item:baseFacet}.
    Assume then that $\Psi_{1} = (\sigma_{1}, \vx_{1}, 0)$ and $\Psi_{2}=(\sigma_{2}, \vx_{2}, 0)$.
    For a vector $v$, let $m_v$ be the smallest number among the absolute values of the entries of $v$.
    Observe that if $v = u$ then $m_v = m_u$. Moreover, we have $m_v = m_{v(\sigma, \vx)}$ for all $\sigma \in S_n$ and $x \in C_2^n$.
    Now, $\kng^{i} \Psi_{1} = (\sigma_{1},\vx_{1},iv_0(\sigma_{1}, \vx_{1}))$ and $\kng^{j} \Psi_{2} = (\sigma_{2},\vx_{2},jv_0(\sigma_{2}, \vx_{2}))$.
    What is more, $m_{iv_0(\sigma_{1}, \vx_{1})} = m_{iv_0} = i$ and $m_{jv_0(\\\sigma_{2},\vx_{2} )} = m_{jv_0} = j$.
    Since $i \neq \pm j$, we see that $m_{iv_0} \neq m_{jv_0}$ and thus $F(\kng^{i} \Psi_{1}) = iv_0(\sigma_{1}, \vx_{1}) \neq jv_0(\sigma_{2}, \vx_{2}) = F(\kng^{j}\Psi_{2})$. We conclude that $\kng^{i} \Psi_{1}$ contains a different face than $\kng^{j} \Psi_{2}$.
    A similar argument proves that  $F(\etab^{i}\Psi_{1}) \neq F(\kng^{j}\Psi_{2})$.
    The remaining cases for Item\nobreakspace \ref {item:NMA} follow from Item\nobreakspace \ref {item:barEta}.

    Items\nobreakspace \ref {item:esAntipoda} and\nobreakspace  \ref {item:esVecina} follow directly from Item\nobreakspace \ref {item:anyFacet} and Item\nobreakspace \ref {item:NMA}.
\end{proof}

Now we are ready to develop our main construction and as a consequence, prove Theorem\nobreakspace \ref {thm:elChido}.
Let $\Fw(\cT)$ denote the set of white flags of $\cT$.
The idea is to build a permutation group $\groupExt$ (acting on the left) on the set $\Fw(\cT) \times \bZ_{3 p}$ for some $p \in \bN$.
We shall prove that the group $\groupExt$ satisfies the relations required to be the rotation group of a rotary extension of $\cT$;
namely those in Corollary\nobreakspace \ref {coro:relsOfExtensionsTaus}.
Then we will impose conditions on $p$ to prove that the potential abstract polytope is actually chiral (and not regular).
Finally, we shall use Lemma\nobreakspace \ref {lem:Lemma10} to prove that the group $\groupExt$ satisfies the intersection property.

We shall denote the elements of $\Fw(\cT) \times \bZ_{3p}$  as pairs $(\Phi,\ell)$ (for $\Phi \in \Fw(\cT)$ and $\ell \in \bZ_{3p}$) keeping in mind that if necessary, $\Phi$ will be represented as a triplet $(\sigma, \vx, t)$, following the labelling of $\Fw(\cT)$ described in Section\nobreakspace \ref {sect:combinatorialToroids}.
Observe that if $l_{1}, l_{2} \in \bZ$ are such that $l_{1}\equiv l_{2} \pmod{3p}$, then $l_{1} \equiv l_{2} \pmod{3}$. Therefore, for $\ell \in \bZ_{3p}$ expressions such as $\ell \equiv 0 \pmod{3}$ are well defined.

Let us establish some notation that we will carry out along the whole section.

Assume that the base flag of $\cT$ is $\Phi_{0}$ and $F_{0}$ is the facet of $\Phi_{0}$.
Let $\kng$ be as in Equation\nobreakspace \textup {(\ref {eq:eta})}.
For $i \in \{\pm 1, \pm 2, \pm 3\}$ let $\Phi_i$ denote the flag $\kng^{i}\Phi_{0}$ and let $F_{i} $ denote the facet of $\Phi_i$.
If $F$ is a facet of $\cT$ we denote by $(F,\ell)$ the set \[\left\{ (\Phi,\ell) : \Phi \text{ is a white flag  and }F \in \Phi \right\}.   \] 
We shall abuse language and call $(F,\ell)$ the \emph{facet} of $(\Phi,\ell)$ whenever $(\Phi,\ell) \in (F,\ell)$.

\begin{defn}\label{defn:hoyo}
Let $X_0$ be the base vertex of $\cT$ and let $X_H$ be the translate of $X_0$ by the vector $(0,-3,0^{n-2})$. Define 
\begin{align*}
H = \{(F,1) : \text{ $F$ is incident to $X_H$} \}.
\end{align*}
\end{defn}

Observe that, by Corollary\nobreakspace \ref {coro:facetsAroundAnyVertex},  if $(F,1) \in H$ then the vector associated with $F$  is of the form $(\epsilon_1,-3 +\epsilon_2, \epsilon_3, \ldots,\epsilon_{n})$ with $\epsilon_i \in \{0,-1\}$. 
From this, we see that if a flag $\Psi$ is in the same facet as $\kng^j\Phi_0$ or $\etab^j\Phi_0$ with $j \in \{0, \pm 1,\pm 2, \pm 3\}$, then the facet of $\Psi$ is not in $H$. This information is relevant for the following definition and observations, but the set $H$ itself will not play a role until Section\nobreakspace \ref {sec:IP}, when we prove that, under certain circumstances, $\groupExt$  has the intersection property.

\begin{defn}\label{defn:baseFlags}
Let $F$ be a facet of $\cT$ and $\ell \in \bZ_{3p}$. 
The \emph{root} $(\Phi,\ell)$ of $(F,\ell)$ is defined as follows:
\begin{enumerate}
    \item \label{item:WhiteBase} If $F = F_i$ for some $i \in \{0,\pm1,\pm2,\pm3\}$, the root of $(F,\ell)$ is $(\Phi_{i}, \ell)$. 
    \item \label{item:BlackBase} If $\ell \not\equiv 1 \pmod{3}$ and $F=F(\etab^{j}\Phi)$ for some white flag $\Phi$ with $F_{0}\in \Phi$ and $j \in \{\pm1,\pm2,\pm3\}$, the root of $\left( F,\ell \right)$  is $(\etab^{j}\Phi,\ell)$.
    \item \label{item:rootsH} If $(F,1) \in H$, then the root of $(F,1)$ is $(\Psi,1)$, with $\Psi$ any white flag incident to the vertex $X_H$.
    \item \label{item:OthersBase} In any other case the root of $(F,\ell)$ is $(\Psi, \ell)$ with $\Psi = (1, \vect{1}, F)$. That is, $\Psi$ is the translate of $\Phi_{0}$ that contains $F$.
    In particular, if  $\ell \equiv 1 \pmod{3}$ and $(F,\ell) \not\in H$, the root of $(F,\ell)$ is $(\Psi, \ell)$ with $\Psi$ a translate of $\Phi_{0}$.
    \end{enumerate}
\end{defn}

\begin{figure}
  \def\svgwidth{1\textwidth }
  \begingroup \makeatletter \providecommand\color[2][]{\errmessage{(Inkscape) Color is used for the text in Inkscape, but the package 'color.sty' is not loaded}\renewcommand\color[2][]{}}\providecommand\transparent[1]{\errmessage{(Inkscape) Transparency is used (non-zero) for the text in Inkscape, but the package 'transparent.sty' is not loaded}\renewcommand\transparent[1]{}}\providecommand\rotatebox[2]{#2}\newcommand*\fsize{\dimexpr\f@size pt\relax}\newcommand*\lineheight[1]{\fontsize{\fsize}{#1\fsize}\selectfont}\ifx\svgwidth\undefined \setlength{\unitlength}{1499.99452125bp}\ifx\svgscale\undefined \relax \else \setlength{\unitlength}{\unitlength * \real{\svgscale}}\fi \else \setlength{\unitlength}{\svgwidth}\fi \global\let\svgwidth\undefined \global\let\svgscale\undefined \makeatother \begin{picture}(1,0.99999988)\lineheight{1}\setlength\tabcolsep{0pt}\put(0,0){\includegraphics[width=\unitlength,page=1]{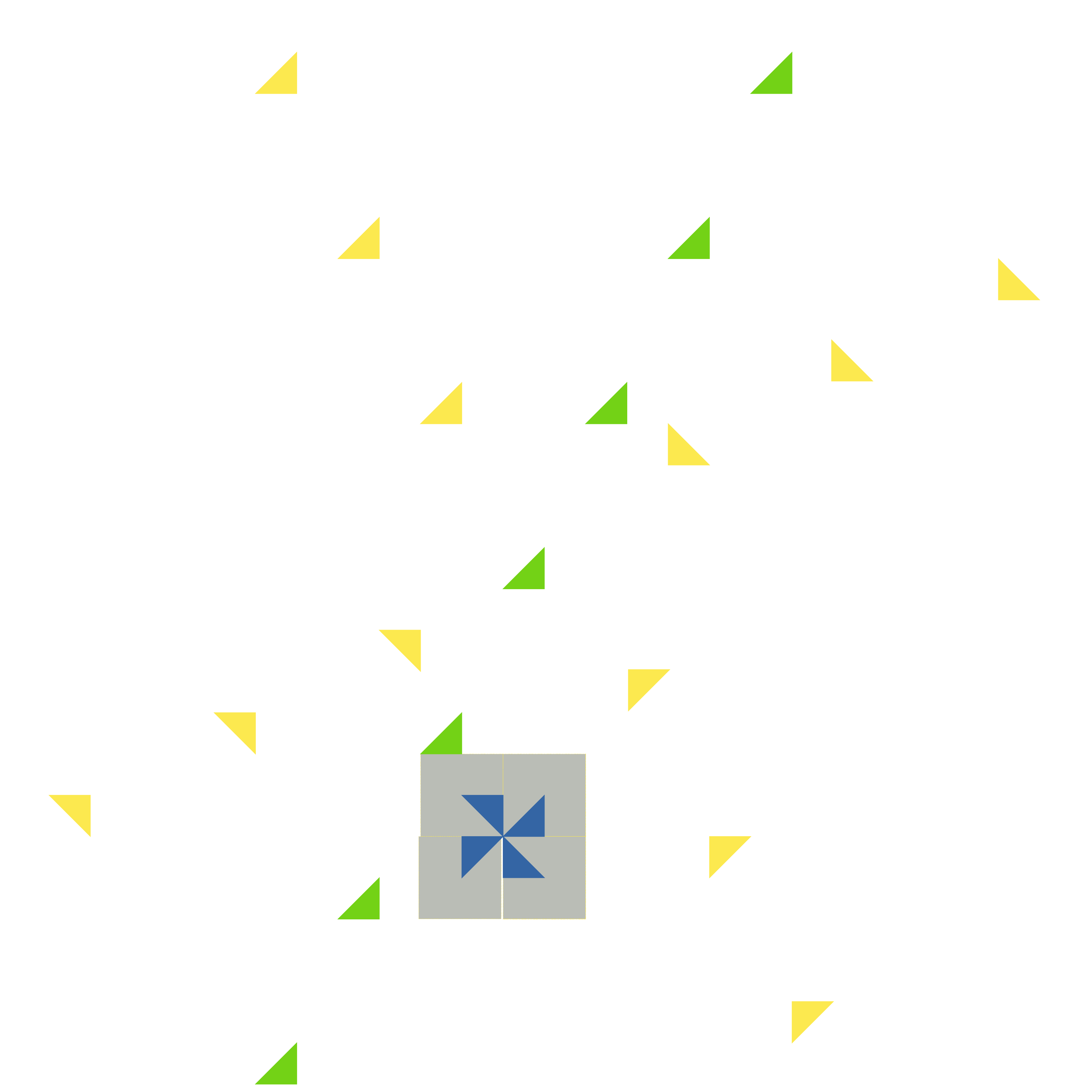}}\put(0.4951715,0.48610869){\color[rgb]{0,0,0}\makebox(0,0)[t]{\lineheight{1.25}\smash{\begin{tabular}[t]{c}$\Phi_0$\end{tabular}}}}\put(0.57406527,0.63433804){\color[rgb]{0,0,0}\makebox(0,0)[t]{\lineheight{1.25}\smash{\begin{tabular}[t]{c}$\Phi_1$\end{tabular}}}}\put(0.64965612,0.78551969){\color[rgb]{0,0,0}\makebox(0,0)[t]{\lineheight{1.25}\smash{\begin{tabular}[t]{c}$\Phi_2$\end{tabular}}}}\put(0.42288363,0.33197473){\color[rgb]{0,0,0}\makebox(0,0)[t]{\lineheight{1.25}\smash{\begin{tabular}[t]{c}$\Phi_{-1}$\end{tabular}}}}\put(0.34729281,0.18079307){\color[rgb]{0,0,0}\makebox(0,0)[t]{\lineheight{1.25}\smash{\begin{tabular}[t]{c}$\Phi_{-2}$\end{tabular}}}}\put(0.72524683,0.93670123){\color[rgb]{0,0,0}\makebox(0,0)[t]{\lineheight{1.25}\smash{\begin{tabular}[t]{c}$\Phi_3$\end{tabular}}}}\put(0.27170196,0.02961142){\color[rgb]{0,0,0}\makebox(0,0)[t]{\lineheight{1.25}\smash{\begin{tabular}[t]{c}$\Phi_{-3}$\end{tabular}}}}\put(0,0){\includegraphics[width=\unitlength,page=2]{chirExtToroids_rootFlags.pdf}}\end{picture}\endgroup    \caption{Roots for $\cT_{\ell}$ and their corresponding $\rho^{(F;\ell)}$ (see Definition\nobreakspace \ref {defn:baseFlags}). The roots of Item\nobreakspace \ref {item:WhiteBase} are shown in green. The flags defined in Item\nobreakspace \ref {item:BlackBase} are shown in yellow and they only apply if $\ell \equiv i \pmod{3}$ for $i \in \left\{ 0,2 \right\}$. The flags described in Item\nobreakspace \ref {item:rootsH} are drawn in blue and they are as shown only if $\ell = 1$. The rest of the roots are translates of $\Phi_{0}$ (
    Item\nobreakspace \ref {item:OthersBase}
  );  some of the associated reflections are drawn in red when they do not coincide with a previously described reflection. }
  \label{fig:rootFlags}
\end{figure}

In Figure\nobreakspace \ref {fig:rootFlags} we represent the choice of the roots for dimension $2$ in the torus $\left\{ 4,4 \right\}_{(13,0)}$.

Note that in most cases the root of the pair $(F,\ell)$ is a translate of $(\Phi_{0},\ell)$.
Strictly speaking, there could be facets $(F, \ell)$ with more than one root.
For example if $n$ is odd, according to Item\nobreakspace \ref {item:WhiteBase} of Definition\nobreakspace \ref {defn:baseFlags} the facet $(F_{1},0)$ has $(\Phi_{1},0)$ as its root; however, according to Item\nobreakspace \ref {item:BlackBase}, the root of $(F_{1},0)$ should be $(r_{0}\hat{\Phi_{1}},0) = (\etab^{-1}r_{0}\hat{\Phi_{0}},0)$.
We let that pass for now but in Remak\nobreakspace \ref {rem:BFfixed} it will become obvious that this potential double definition of roots is not relevant for our purposes.

Now, consider a facet $F$ of $\cT$ and an integer $\ell \in \bZ_{3p}$, and suppose $(\Phi_F,\ell)$ is a root of $(F,\ell)$. For $i \in \{0,1\}$, we define $\rho_i^{(F,\ell)}$  as the automorphism of $F$ mapping $\Phi_F$ to $r_i\Phi_F$. As we will see shortly, this permutation does not depend on the choice of the root of $(F,\ell)$ (if there happens to be more than one). Further, $\rho_i^{(F,\ell)}$ permutes the flags of $F$ and maps white flags to black flags and vice versa. For every $\ell \in \bZ_{3p}$, we can define a permutation $\rho^\ell$ of $\cF(\cT)$
 by 
 \begin{equation}\label{eq:rr} \rho^\ell \Phi = 
 \begin{cases}
 \Phi \rho_0^{(F,\ell)},   
 &\text{if $ (F,\ell) \notin H$}\\
 \Phi \rho_1^{(F,\ell)},
 &\text{if $ (F,\ell) \in H$}
 \end{cases}
 \end{equation}
 where $F$ is the facet of $\Phi$.
A fairly obvious observation is that for every $\ell$, $\rr$ is an involution.
Moreover, we have the following results.

\begin{rem}\label{rem:rchiquitas}
  Let $0 \leq i \leq n-1$ and $\ell \in \bZ_{3p}$, $\Phi$ a flag of $\cT$ and $F$ the facet of $\Phi$.
  Since the element $r_{i}$ preserves the facet of $\Phi$ and both $\rho_{0}^{(F,\ell)}$ and $\rho_{1}^{(F,\ell)}$ are automorphisms of $F$, then 
  \begin{align*}
    r_{i}\left( \Phi \rho_{0}^{(F,\ell)} \right) &=\left(r_{i}  \Phi \right) \rho_{0}^{(F,\ell)};  \text{ and} \\ 
    r_{i}\left( \Phi \rho_{1}^{(F,\ell)} \right) &=\left(r_{i}  \Phi \right) \rho_{1}^{(F,\ell)}. 
  \end{align*}
  It follows that $\rr$ and $r_{i}$ commute.
\end{rem}

\begin{rem}\label{rem:BFfixed_}
   Let $F$ be a facet of $\cT$ and $\ell \in \bZ_{3p} $ such that $(F,\ell) \not\in H$.
   Let $\Phi_{F}$ be the flag containing $F$ such that $(\Phi_{F}, \ell)$ is a root of $(F,\ell)$.
   If $\Phi \in \{\Phi_{F}, \hat{\Phi_F}, r_{0} \hat{\Phi_F}  \}$, then $\rr\Phi = r_{0}\Phi$, or equivalently,
   \begin{equation}\label{eq:BFfixed}
   \rr r_{0} (\Phi) = \Phi.
   \end{equation}
 \end{rem}

 Following the notation in the remark above, all the flags $\Phi_{F}$, $ \hat{\Phi_{F}}$ and $r_{0} \hat{\Phi_{F}}$ define the same $\rho^{(F,\ell)}$ and thus $\rr$ is well defined regardless of the potential existence of two base flags for a facet $F$.
 Moreover, the only property of those base flags that will be relevant to us, is the one described in Remak\nobreakspace \ref {rem:BFfixed}.

We are now ready to give explicit  a definition for generators of the group $\groupExt = \langle \te_1, \dots, \te_{n+1} \rangle$. For $ 1\leq i \leq n$ define.
\begin{equation} \label{eq:smallGens}
\te_{i} \left(\Phi, \ell \right) = \left(r_{0}r_{i}\Phi,\ell \right).
\end{equation}The permutation $\te_{n+1}$ is given by 
\begin{equation}\label{eq:lastGen}
\te_{n+1} \left(\Phi, \ell \right) =
  \begin{cases}
	  \left(\rr r_{0}\Phi,\ell+1 \right), 
            &\text{if } F_{1} \in \Phi\text{ and } \ell \equiv 0 \pmod{3};\\ 
            &\text{or } F_{3} \in \Phi \text{ and } \ell \equiv 1 \pmod{3}; \\
            &\text{or } F_{0} \in \Phi \text{ and } \ell \equiv 2 \pmod{3}; \\
	  \left(\rr r_{0}\Phi,\ell-1 \right), 
            &\text{if } F_{0} \in \Phi\text{ and } \ell \equiv 0 \pmod{3};\\ 
            &\text{or } F_{1} \in \Phi \text{ and } \ell \equiv 1 \pmod{3}; \\
            &\text{or } F_{3} \in \Phi \text{ and } \ell \equiv 2 \pmod{3}; \\
	  \left(\rr r_{0}\Phi,\ell \right), &\text{otherwise.}
 \end{cases}
\end{equation}

We can rewrite Remak\nobreakspace \ref {rem:BFfixed_} as follows.

\begin{remark}\label{rem:BFfixed}
 Let $F$ be a a facet of $\cT$ and $\ell \in \bZ_{3p} $ such that $(F,\ell) \not\in H$, then for all  
 for all $\Phi \in \{\Phi_{F}, \hat{\Phi_F}, r_{0} \hat{\Phi_F}  \}$ we have
 \begin{equation}
  \te_{n+1} (\Phi,\ell) = (\Phi,\ell')
\end{equation}
where 
\[ 
  \ell' =
  \begin{cases}
	  \ell+1, 
            &\text{if } F= F_{1} \text{ and } \ell \equiv 0 \pmod{3};\\ 
            &\text{or } F=F_{3}\text{ and } \ell \equiv 1 \pmod{3}; \\
            &\text{or } F= F_{0} \text{ and } \ell \equiv 2 \pmod{3}; \\
	  \ell-1 , 
            &\text{if } F = F_{0} \text{ and } \ell \equiv 0 \pmod{3};\\ 
            &\text{or } F = F_{1}  \text{ and } \ell \equiv 1 \pmod{3}; \\
            &\text{or } F = F_{3} \text{ and } \ell \equiv 2 \pmod{3}; \\
	  \ell, &\text{otherwise.}
 \end{cases}
\] 
\end{remark}

Now, note that for a fixed $\ell \in \bZ_{3p}$, the group $\groupExt_{n+1} := \langle \te_{1}, \dots, \te_{n} \rangle$ acts on  $\cT_\ell$ just as the group $\langle r_{0}r_{1}, r_{0} r_{2}, \dots, r_{0}r_{n} \rangle \cong \monp(\cT)$ acts on the set $\Fw(\cT)$.
We abuse  notation by considering the elements of $\monw(\cT)$ as elements of $\groupExt$. 
In particular, if $\kng$ denotes the element defined in Equation\nobreakspace \textup {(\ref {eq:eta})}, then we may think that $\kng \in \groupExt$.
Moreover, $\groupExt_{n+1} \cong \monw(\cT) \cong \autp(\cT)$ (see Part\nobreakspace \ref {part:isomorphismChir} of Proposition\nobreakspace \ref {prop:autpAsMonPChir}).
This implies that the (potential) polytope defined by $\groupExt$ has facets isomorphic to $\cT$.

The group elements $\te_{1}, \dots \te_{n+1}$ will play the role of the elements $\tau_{1}, \dots, \tau_{n+1}$ of Corollary\nobreakspace \ref {coro:relsOfExtensionsTaus} (note the shift of indices), meaning that in order to guarantee that the group $\groupExt$ is the automorphism group of a chiral extension of $\cT$, we need to prove that the generators $\te_{1}, \dots, \te_{n+1}$ satisfy the relations in Equation\nobreakspace \textup {(\ref {eq:relsOfExtensionsTaus})}.

\begin{prop} \label{prop:chiralExtRelations}
 With the notation given above, the group elements $\langle \te_{1}, \dots, \te_{n+1} \rangle$ satisfy the relations of Corollary\nobreakspace \ref {coro:relsOfExtensionsTaus}.
\end{prop}
\begin{proof}
 Translating the notation of Corollary\nobreakspace \ref {coro:relsOfExtensionsTaus}, we only need to prove that
 \[\begin{aligned}
  \te_{n+1}^{2} &= 1, && \\
 (\te_{i}^{-1} \te_{n+1})^{2} &= 1  &&\text{for all } 1 \leq i \leq n-1.
\end{aligned}\]

First note that an immediate consequence of Remak\nobreakspace \ref {rem:rchiquitas} is that $\rr r_{0}$ is a transposition.
Assume that $\Phi$ is a white flag with facet $F_{1}$ and $\ell \in \bZ_{3p}$ is such that $\ell \equiv 0 \pmod{3}$.
Then \[\te_{n+1}^{2}(\Phi, \ell) = \te_{n+1}(\rr r_{0} \Phi, \ell+1),\] but $\rr r_{0}$ maps flags containing $F_{1}$ to flags containing $F_{1}$ and $\ell+1 \equiv 1 \pmod{3}$, hence \[\te_{n+1}(\rr r_{0} \Phi, \ell+1) = \left( \rr[\ell +1] r_{0} \rr r_{0} \Phi, \ell\right).\]
Finally observe that $\rho_{0}^{(F_{1},\ell)} = \rho_{0}^{(F_{1},\ell+1)}$; which implies that $\rr[\ell +1] r_{0} \rr r_{0} \Phi = \Phi$.
With a similar argument, we can prove that if $\Phi$ and $\ell$ are any of the other possible pairs such that the second coordinate of $\te_{n+1}(\Phi,\ell)$ is not $\ell$, then $\te_{n+1}^2 (\Phi,\ell) = \left( \Phi, \ell \right)$.
If $\te_{n+1}(\Phi,\ell) = (\rr r_0 \Phi,\ell)$ then clearly $\te^{2}_{n+1}(\Phi,\ell) = (\Phi,\ell)$.

To prove the second part, first note that if $\Phi$ is a flag and $1 \leq i \leq n-1$, then $\Phi$ and $r_{i} \rr \Phi$ contain the same facet. 
Again, we will only prove the case when $F_{1} \in \Phi$ and $\ell \equiv 0 \pmod{3}$; the remaining cases are analogous. We have
\begin{align*}
 (\te_{i}^{-1} \te_{n+1})^{2}\left(\Phi, \ell\right) &= \te_{i}^{-1} \te_{n+1} \left(r_{i}r_{0} \rr r_{0}\Phi, \ell+1\right)  \\
 &=  \te_{i}^{-1} \te_{n+1} \left(r_{i}\rr \Phi, \ell+1\right) \\
 &= \te_{i}^{-1}(\rr[\ell + 1] r_{0} r_{i} \rr \Phi, \ell) \\
 &= (r_{i} r_{0} \rr[\ell + 1] r_{0} r_{i} \rr \Phi, \ell) \\
 &= (\Phi, \ell).
\end{align*}
The last equality follows from the facts that $r_{0}$ and $r_{i}$ commute with $\rr$ (Remak\nobreakspace \ref {rem:rchiquitas}) and $\rr$ and $\rr[\ell+1]$ induce the same permutation on the flags containing $F_{1}$.
\end{proof}

Next we will prove that under certain conditions, there is no group automorphism as described in Corollary\nobreakspace \ref {cor:gpAutAltGens}. 
More precisely,

\begin{prop}\label{prop:noGroupAut}
  Let $\groupExt$ be the group generated by the permutations $ \te_{1}, \dots, \te_{n+1} $ defined in Equations\nobreakspace \textup {(\ref {eq:smallGens})} and\nobreakspace  \textup {(\ref {eq:lastGen})} acting on  the set $\Fw(\cT) \times \bZ_{3p}$, for a prime $p > 3|\Fw(\cT)|$.
  Let
  \begin{equation}
    \label{eq:mu}
    \begin{aligned}
      \mu &= \te_{n+1}\kng^{-3}\te_{n+1} \kng^{2} \te_{n+1} \kng &&\text{and} \\
      \mub&= r_{0}\kng^{-3} r_{0}\te_{n+1} r_{0}\kng^{2} r_{0} \te_{n+1} r_{0}\kng r_{0}.
    \end{aligned}
  \end{equation}

  Then, there is no group automorphism $\alpha: \groupExt \to \groupExt$ such that \[\alpha(\mu) = \mub.\]
  In particular, there is no automorphism $\alpha:\groupExt \to \groupExt $ satisfying
  \[\begin{aligned}
	    \alpha(\te_{1}) &= \te^{-1}_{1},&& \\
    \alpha(\te_{i}) &= \te_{i} &&\text{for all } 2 \leq i \leq n+1.
\end{aligned}\]
\end{prop}

The idea behind the proof goes as follows: the hypothetical automorphism $\alpha$ acts on $\left\langle    \te_{1}, \dots, \te_{n} \right\rangle \cong \monw(\cT) $ as conjugation by $r_{0}$.
Since $\alpha$ fixes $\te_{n+1}$, it should map $\mu$ to $\mub$.
To prove that this is impossible we shall prove that $\langle \mu \rangle$ induces an orbit of length $p$ (Lemma\nobreakspace \ref {lem:orbitMu}) and that the size of the orbits $\langle \mub \rangle$ is bounded (Lemma\nobreakspace \ref {lem:onlyThreeCopies}).

By choosing $p$ equal to a large prime we would have proved that $p$ divides the order of $\mu$ but does not divide the order of $\mub$, which in turn implies that $\alpha$ cannot map $\mu$ to $\mub$.

\begin{figure}\centering
\def\svgwidth{\textwidth}
\begin{scriptsize}
	\begingroup \makeatletter \providecommand\color[2][]{\errmessage{(Inkscape) Color is used for the text in Inkscape, but the package 'color.sty' is not loaded}\renewcommand\color[2][]{}}\providecommand\transparent[1]{\errmessage{(Inkscape) Transparency is used (non-zero) for the text in Inkscape, but the package 'transparent.sty' is not loaded}\renewcommand\transparent[1]{}}\providecommand\rotatebox[2]{#2}\newcommand*\fsize{\dimexpr\f@size pt\relax}\newcommand*\lineheight[1]{\fontsize{\fsize}{#1\fsize}\selectfont}\ifx\svgwidth\undefined \setlength{\unitlength}{2633.85278897bp}\ifx\svgscale\undefined \relax \else \setlength{\unitlength}{\unitlength * \real{\svgscale}}\fi \else \setlength{\unitlength}{\svgwidth}\fi \global\let\svgwidth\undefined \global\let\svgscale\undefined \makeatother \begin{picture}(1,0.26104334)\lineheight{1}\setlength\tabcolsep{0pt}\put(0,0){\includegraphics[width=\unitlength,page=1]{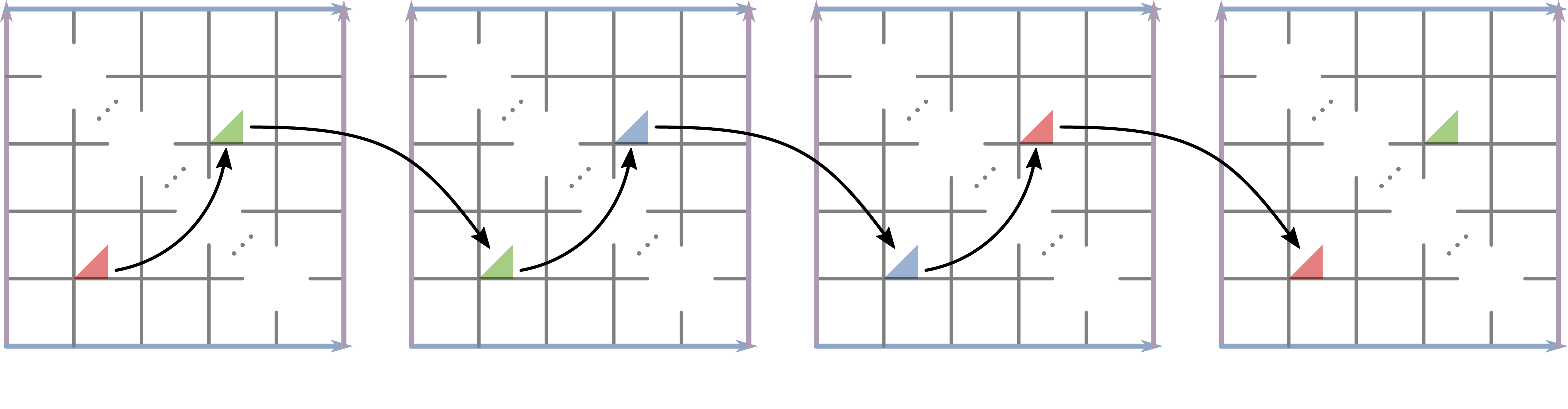}}\put(0.13324731,0.10998721){\color[rgb]{0,0,0}\makebox(0,0)[lt]{\lineheight{1.25}\smash{\begin{tabular}[t]{l}$\kng$\end{tabular}}}}\put(0.24087085,0.18532369){\color[rgb]{0,0,0}\makebox(0,0)[t]{\lineheight{1.25}\smash{\begin{tabular}[t]{c}$\te_{n+1}$\end{tabular}}}}\put(0.39154381,0.10998721){\color[rgb]{0,0,0}\makebox(0,0)[lt]{\lineheight{1.25}\smash{\begin{tabular}[t]{l}$\kng^{2}$\end{tabular}}}}\put(0.64984031,0.10998721){\color[rgb]{0,0,0}\makebox(0,0)[lt]{\lineheight{1.25}\smash{\begin{tabular}[t]{l}$\kng^{-3}$\end{tabular}}}}\put(0.1117226,0.00236367){\color[rgb]{0,0,0}\makebox(0,0)[t]{\lineheight{1.25}\smash{\begin{tabular}[t]{c}$\cT_{\ell}$\end{tabular}}}}\put(0.37001995,0.00236368){\color[rgb]{0,0,0}\makebox(0,0)[t]{\lineheight{1.25}\smash{\begin{tabular}[t]{c}$\cT_{\ell+1}$\end{tabular}}}}\put(0.62831639,0.00236368){\color[rgb]{0,0,0}\makebox(0,0)[t]{\lineheight{1.25}\smash{\begin{tabular}[t]{c}$\cT_{\ell+2}$\end{tabular}}}}\put(0.88661289,0.00236368){\color[rgb]{0,0,0}\makebox(0,0)[t]{\lineheight{1.25}\smash{\begin{tabular}[t]{c}$\cT_{\ell+3}$\end{tabular}}}}\put(0.49916735,0.18532369){\color[rgb]{0,0,0}\makebox(0,0)[t]{\lineheight{1.25}\smash{\begin{tabular}[t]{c}$\te_{n+1}$\end{tabular}}}}\put(0.75746385,0.18532369){\color[rgb]{0,0,0}\makebox(0,0)[t]{\lineheight{1.25}\smash{\begin{tabular}[t]{c}$\te_{n+1}$\end{tabular}}}}\put(0.06329201,0.10460604){\makebox(0,0)[t]{\lineheight{1.25}\smash{\begin{tabular}[t]{c}$\Phi_0$\end{tabular}}}}\put(0.14939084,0.19070487){\makebox(0,0)[t]{\lineheight{1.25}\smash{\begin{tabular}[t]{c}$\Phi_1$\end{tabular}}}}\put(0.32158851,0.10460604){\makebox(0,0)[t]{\lineheight{1.25}\smash{\begin{tabular}[t]{c}$\Phi_1$\end{tabular}}}}\put(0.40768734,0.19070487){\makebox(0,0)[t]{\lineheight{1.25}\smash{\begin{tabular}[t]{c}$\Phi_3$\end{tabular}}}}\put(0.57988501,0.10460604){\makebox(0,0)[t]{\lineheight{1.25}\smash{\begin{tabular}[t]{c}$\Phi_3$\end{tabular}}}}\put(0.66598384,0.19070487){\makebox(0,0)[t]{\lineheight{1.25}\smash{\begin{tabular}[t]{c}$\Phi_0$\end{tabular}}}}\put(0.83818151,0.10460604){\makebox(0,0)[t]{\lineheight{1.25}\smash{\begin{tabular}[t]{c}$\Phi_0$\end{tabular}}}}\end{picture}\endgroup  \end{scriptsize}
\caption{The pair $(\Phi_{0}, \ell)$ under the action of $\mu$ \label{fig:jumps}}
\end{figure}
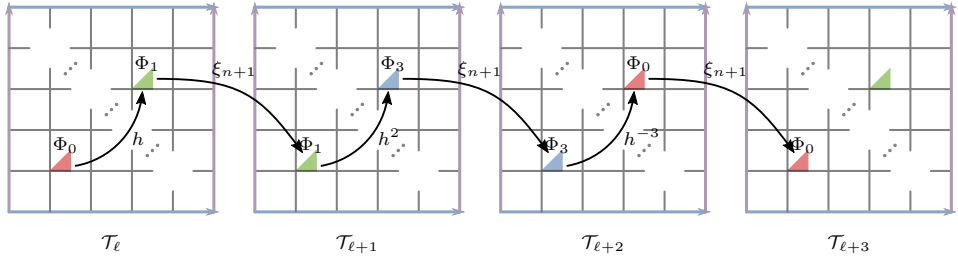

\begin{lem}\label{lem:orbitMu}
 Let $\Phi_{0}$ denote the base flag of $\cT$. Let $\groupExt = \left\langle \te_{1}, \dots, \te_{n+1} \right\rangle $ be the permutation group on the set $\Fw(\cT) \times \bZ_{3p}$ defined by Equations\nobreakspace \textup {(\ref {eq:smallGens})} and\nobreakspace  \textup {(\ref {eq:lastGen})}.
 Let $\mu$ be as in Equation\nobreakspace \textup {(\ref {eq:mu})}.
 Then the orbit of $(\Phi_{0}, 0)$ under $\langle \mu \rangle$ has $p$ elements.
\end{lem}
\begin{proof} 
    As before, for $i\in\{\pm1, \pm 2, \pm 3\}$ let $\Phi_i$ denote the flag $\kng^{i}\Phi_0$ and let $F_i$ denote the facet of $\Phi_i$.
    Recall that for every $j  \in \bZ_{3p}$,  $\te_{n+1}(\Phi_{i},j) = (\Phi_{i},j')$ with $j' \in \{j-1, j , j+1\}$ (see Remak\nobreakspace \ref {rem:BFfixed}).
    Let $\ell \in \bZ_{3p}$ be such that $\ell \equiv 0 \pmod{3}$.
  Consider then the following computation:
\begin{subequations}
 \begin{align}
  \mu(\Phi_{0},\ell)
  &= \te_{n+1} \kng^{-3} \te_{n+1} \kng^{2} \te_{n+1} \kng (\Phi_{0}, \ell) \nonumber \\
  &= \te_{n+1} \kng^{-3} \te_{n+1} \kng^{2} \te_{n+1} (\Phi_{1}, \ell)  \nonumber \\
  &= \te_{n+1} \kng^{-3} \te_{n+1} \kng^{2}  (\Phi_{1}, \ell+1)  \label{eq:orbitMu_1}  \\
&= \te_{n+1} \kng^{-3} \te_{n+1}  (\Phi_{3}, \ell+1)  \nonumber  \\
  &= \te_{n+1} \kng^{-3}   (\Phi_{3}, \ell+2)   \label{eq:orbitMu_2} \\
&= \te_{n+1}   ( \Phi_{0}, \ell+2)  \nonumber \\
  &= (\Phi_{0}, \ell+3),  \label{eq:orbitMu_3} 
\end{align}
\end{subequations}
were in \eqref{eq:orbitMu_1}, \eqref{eq:orbitMu_2} and \eqref{eq:orbitMu_3}, we used that $F_{i} \in \Phi_{i}$ for $i = 1, 3$ and $0$, respectively (see Figure\nobreakspace \ref {fig:jumps}
 and Equation\nobreakspace \textup {(\ref {eq:lastGen})}).
 It follows that $\mu^{j}(\Phi_0, 0) = (\Phi_0, 3 j)$, which implies that the length of the orbit of $(\Phi_0, 0)$ under $\langle \mu \rangle$ is $p$.
\end{proof}

Let us now explore the possibilities for the orbits under $\left\langle \mub \right\rangle $. Informally, we aim to show that the orbit of a flag under  $\left\langle \mub \right\rangle$ is either very small, or is confined within three `copies' of $\mathcal{T}$. This is stated formally in Lemma\nobreakspace \ref {lem:onlyThreeCopies} below, which will be proved at the end of this section.

\begin{lem}\label{lem:onlyThreeCopies}
  Let $(\Phi,\ell) \in \Fw(\cT) \times \bZ_{3p}$ and let $t \in \bZ_{3p}$ such that $\ell \in \left\{ 3t, 3t+1, 3t+2 \right\} $ then one of the following hold:
  \begin{enumerate}
    \item The orbit $\left\langle \mub \right\rangle (\Phi,\ell) $ is contained in $\cT_{3t} \cup \cT_{3t+1} \cup \cT_{3t+2}$; or
    \item $\left| \left\langle \mub \right\rangle (\Phi,\ell)  \right| \leq 3 $.
  \end{enumerate}
\end{lem}

To prove Lemma\nobreakspace \ref {lem:onlyThreeCopies}, we will need a series of auxiliary results. The idea of the proof is as follows. Assume that the orbit of $(\Phi,\ell)$ under $\left\langle \mub \right\rangle $ is not contained in $\cT_{3t} \cup \cT_{3t+1} \cup \cT_{3t+2}$ and let $j$ be such that $\mub^{j}(\Phi, \ell) = (\Psi,\ell') \in \cT_{3t} \cup \cT_{3t+1} \cup \cT_{3t+2}$ but $\mub^{j+1}(\Phi,\ell) = \mub(\Psi,\ell') \not \in \cT_{3t} \cup \cT_{3t+1} \cup \cT_{3t+2}$.
Consider the following pairs
\[\begin{aligned}
	(\Psi_{1}, \ell'_{1}) &= \etab(\Psi,\ell') \\
	(\Psi_{2}, \ell'_{2}) &= \etab^{2}\te_{n+1}\etab(\Psi,\ell') \\
	(\Psi_{3}, \ell'_{3}) &= \etab^{-3}\te_{n+1}\etab^{2}\te_{n+1}\etab(\Psi,\ell') \\
\end{aligned}\]
Since $\mub(\Psi,\ell') \not \in \cT_{3t} \cup \cT_{3t+1} \cup \cT_{3t+2}$, one of the pairs $(\Psi_{i}, \ell'_{i})$ must be such that the second coordinate of $\te_{n+1}(\Psi_{i}, \ell'_{i})$ is equal to $3t-1$ or $3t+3$. Note that this happens if and only if  the pair $(\Psi_{i}, \ell'_{i})$ satisfies the following two conditions:
\begin{enumerate}
    \item $\ell'_{i}\in\{3t,3t+2\}$,
    \item $F_{0} \in \Psi_{i}$.
\end{enumerate}
The orbit of $(\Phi, \ell)$  will depend on which of the three pairs $(\Psi_{i}, \ell'_{i})$ satisfies the above conditions. 
In the following pages, we will show that $(\Phi,\ell)$ 
has an orbit of size $2$ (Lemmas\nobreakspace  \ref {lem:chapulines1} to\nobreakspace  \ref {lem:chapulines3} ) unless $\Phi$ satisfies very specific conditions (Lemma\nobreakspace \ref {lem:chapulines0}), in which case $\left|   \left\langle \mub \right\rangle (\Phi, \ell) \right| = 3$. 
This will imply that  $\left|   \left\langle \mub \right\rangle (\Psi, \ell') \right|  =  \left|   \left\langle \mub \right\rangle (\Phi, \ell) \right| \leq 3$, completing the proof.

\begin{lemma}\label{lem:chapulines0}

Let $\ell \in \bZ_{3p}$ such that $\ell \equiv 0 \pmod{3}$ and $\Phi \in \{r_{0}\hat{\Phi_{1}}, r_{0}\hat{\Phi_{3}} \}$. Observe that $\Phi$ is white if and only if $n$ is odd.
In this situation, the orbits of $(\Phi,\ell)$ under $\left\langle \mub \right\rangle $ are
\[
\begin{gathered}
\left\{ (r_{0}\hat{\Phi_{1}},\ell), (r_{0}\hat{\Phi_{1}},\ell-1), (r_{0}\hat{\Phi_{1}},\ell+1) \right\}, \\
\left\{ (r_{0}\hat{\Phi_{3}},\ell), (r_{0}\hat{\Phi_{3}},\ell-2), (r_{0}\hat{\Phi_{3}},\ell-1) \right\}.
\end{gathered}
\]
In particular, this implies that if  $\ell \equiv i \pmod{3}$ for some $i \in \left\{ 0,2 \right\} $ and $\Phi \in \{r_{0}\hat{\Phi_{1}},r_{0}\hat{\Phi_{3}}\}  $, then
\[
	\left| \left\langle \mub \right\rangle \left( \Phi \right) \right| = 3.
\]
\end{lemma}

\begin{proof}
  Observe that a direct consequence of Equations\nobreakspace \textup {(\ref {eq:antipoda})} and\nobreakspace  \textup {(\ref {eq:etaBar})} is that $\etab^{j} r_{0}\hat{\Phi_{i}} = r_{0}\hat{\Phi_{i-j}}$.
  This can be used to see that the first coordinates of $(\Phi,\ell)$, $(\etab \Phi,\ell)$ and $(\etab^3 \Phi,\ell)$ are all fixed by $\te_{n+1}$, for every $\ell \in \bZ_{3p}$ (Remak\nobreakspace \ref {rem:BFfixed}), which implies that the first coordinate of $\mu(\Phi,\ell)$ is $\Phi$.
  It only remains to track the second coordinate of $\mu(\Phi,\ell)$. 
  We shall only do it for $\Phi = r_{0} \hat{\Phi_1}$, the other case is analogous.
  The second coordinate of $\mu(\Phi,\ell)$ must be $\ell -1$ because $F_{0} \in \etab \Phi$  but the facets of $\etab^{3} \Phi$ and $\Phi$ are not $F_3$ or $F_1$ (note that $\ell-1 \equiv 2 \pmod{3}$).
  A similar argument proves that $\mub^{2}(\Phi, \ell) = (\Phi, \ell + 1)$ and that $\mub^{3}(\Phi,\ell) = (\Phi,\ell)$.
\end{proof}

\begin{lemma}
\label{lem:chapulines1}
Let $(\Phi,\ell) \neq (r_0 \hat{\Phi_1},\ell)$ and let $(\Psi,\ell') =  \etab (\Phi,\ell)$.
Suppose $F_0 \in \Psi$ and $\ell' \equiv i \pmod 3$ for some integer $i\in \{0,2\}$.
Then the following hold:
\begin{enumerate}
    \item \label{item:chapulines11} $(\Phi,\ell) = (\etab^{-1} \Psi, \ell')$, in particular $\ell' = \ell$;
    \item \label{item:chapulines12} if $i = 2$ then $\mub(\Phi,\ell) = (\kng^{-1} \rr r_0 \etab^1 \Phi,\ell+1)$;
    \item \label{item:chapulines13} if $i = 0$ then $\mub(\Phi,\ell) = (\kng^{-1} \rr r_0 \etab^1 \Phi,\ell-1)$;
    \item \label{item:chapulines14}$ |\left\langle \mub \right\rangle (\Phi,\ell)| = 2$.
\end{enumerate}
\end{lemma}

\begin{proof}
  Item\nobreakspace \ref {item:chapulines11} is obvious.
  Assume that $i = 2 $
  and consider the following computation:
  \begin{subequations}\label{eq:chapulines1_Ida}
  \begin{align}
    \mub \left(\Phi, \ell \right)
     &= \te_{n+1} \etab^{-3} \te_{n+1} \etab^{2} \te_{n+1} \etab \left( \Phi, \ell \right) \nonumber \\
     &= \te_{n+1} \etab^{-3} \te_{n+1} \etab^{2} \te_{n+1}  \left(\etab \Phi, \ell \right) \nonumber \\
     &= \te_{n+1} \etab^{-3} \te_{n+1} \etab^{2} \left(\rr r_{0}\etab \Phi, \ell+1  \right)  \label{eq:chapulines1_a}\\
     &= \te_{n+1} \etab^{-3} \te_{n+1}  \left( \etab^{2} \rr r_{0}\etab \Phi, \ell+1  \right) \nonumber \\
&= \te_{n+1} \etab^{-3}\left( \etab^{2} \rr r_{0}\etab \Phi, \ell+1  \right) \label{eq:chapulines1_b}\\
     &= \te_{n+1} \left(\etab^{-1} \rr r_{0}\etab \Phi, \ell+1  \right) \nonumber\\
&= \left(\etab^{-1} \rr r_{0}\etab \Phi, \ell+1  \right), \label{eq:chapulines1_c}
  \end{align}
  \end{subequations}
  where \eqref{eq:chapulines1_a} holds because $F_{0}$ is the facet of $\etab \Phi$.
  The facet of $\etab^{2} (\rr r_{0} \etab \Phi)$ cannot be $F_{0}$ or $F_{1}$ (see Item\nobreakspace \ref {item:NMA} of Proposition\nobreakspace \ref {prop:etaIsGood}) and $\etab^{2} (\rr r_{0} \etab \Phi,\ell+1)$ is the root of its facet; hence \eqref{eq:chapulines1_b} holds.
  Finally observe that $F_{0} \not \in \etab^{-1}(\rr r_{0}\etab \Phi)$, because $F_{0} \in \rr r_{0}\etab \Phi$ and if $F_{1} \in \etab^{-1}(\rr r_{0}\etab \Phi)$, then $\Phi$ must be
  $r_{0}\hat{\Phi_{1}}$, which was excluded by hypothesis.
  Therefore, Equation\nobreakspace \textup {(\ref {eq:chapulines1_c})} holds.
  This proves Item\nobreakspace \ref {item:chapulines12}.

  A completely analogous computation to that in Equation\nobreakspace \textup {(\ref {eq:chapulines1_Ida})} can be used to prove Item\nobreakspace \ref {item:chapulines13}. That is, that if $i=0$, then $\mub(\Phi, \ell) = (\etab^{-1} \rr r_{0} \etab \Psi, \ell-1)$.
  Keep in mind that in this case we shall observe that $F_3$ is not the facet of $ \etab^{2} \rr r_0 \etab \Phi $ (to prove Equation\nobreakspace \textup {(\ref {eq:chapulines1_b})}) or of $ \etab^{-1} \rr r_0 \etab \Phi $ (to prove Equation\nobreakspace \textup {(\ref {eq:chapulines1_c})}), but both cases follow from Item\nobreakspace \ref {item:NMA} of Proposition\nobreakspace \ref {prop:etaIsGood}.

  Let $\Phi' = \etab^{-1} \rr r_{0} \etab \Phi$ and observe that if $\ell \equiv 2 \pmod{3}$, then $(\Phi',\ell+1)$ satisfies the hypotheses of Item\nobreakspace \ref {item:chapulines13}.
  Indeed, $\Phi' \neq r_{0} \hat{\Phi_{1}}$, otherwise $\rr r_{0} \Phi =  r_{0} \hat{\Phi_{0}}$ which in turn implies that $\Phi = r_{0} \hat{\Phi_{1}}$.
  Clearly $\ell+1 \equiv 0 \pmod{3}$ and $F_{0} \in \kng \Phi' = \rr r_{0} \etab \Phi$.
  A similar argument can be used to show that if $\ell \equiv 0 \pmod{3}$, then $(\Phi',\ell-1)$ satisfies the hypotheses of Item\nobreakspace \ref {item:chapulines12}.

  The proof is complete by observing that
  \[\begin{aligned}
  \etab^{-1} \rr[\ell+1] r_{0} \etab
  (\etab^{-1} \rr r_{0} \etab \Phi) &= \Phi &\text{and} \\
  \etab^{-1} \rr[\ell-1] r_{0} \etab (\etab^{-1} \rr r_{0} \etab \Phi) &= \Phi,
    \end{aligned}
  \]
  where in both cases we have used that $\rr[\ell-1]$, $ \rr$ and $\rr[\ell+1]$ act the same way on the flags containing $F_0$.
\end{proof}

\begin{lemma}\label{lem:chapulines2}
Let $(\Phi,\ell) \neq (r_0 \hat{\Phi_3},\ell)$ and let $(\Psi,\ell') =  \etab^{2} \te_{n+1} \etab (\Phi,\ell)$.
Suppose $F_0 \in \Psi$ and $\ell' \equiv i \pmod 3$ for some integer $i\in \{0,2\}$. 
Then the following hold:
\begin{enumerate}
    \item \label{item:chapulines21} $(\Psi,\ell') = (\etab^{3}\Phi, \ell)$;
\item \label{item:chapulines22} if $i = 2$ then $\mub(\Phi,\ell) = (\etab^{-3} \rr r_0 \etab^3 \Phi,\ell+1)$;
    \item \label{item:chapulines23} if $i = 0$ then $\mub(\Phi,\ell) = (\etab^{-3} \rr r_0 \etab^3 \Phi,\ell-1)$;
    \item \label{item:chapulines24}$ |\left\langle \mub \right\rangle (\Phi,\ell)| = 2$. 
\end{enumerate}

\end{lemma}

\begin{proof}
 Let $(\Phi,\ell)$ and $(\Psi, \ell')$ as described above.
 Since $F_{0} \in \Psi$, then clearly the facet of $\etab^{-2}\Psi$ is not $F_{0}$, $ F_{1} $ or $F_{3}$, which implies
 \[\etab\left( \Phi,\ell \right) = \te_{n+1} \etab^{-2} (\Psi, \ell') = (\etab^{-2} \Psi, \ell' ),\]
 where the last equality follows from the fact that $\etab^{-2}(\Psi, \ell' )$ is the root of its facet.
 This proves that $\etab (\Phi,\ell) = (\etab^{-2 } \Psi,\ell')$, or equivalently $(\Psi,\ell) = (\etab^{3 } \Phi,\ell')$.

  Now assume that $\ell \equiv 2 \pmod{3}$.
  Let us now compute $\mub\left( \Phi, \ell \right)$:
  \begin{subequations} \label{eq:chapulines2}
    \begin{align}
      \mub\left( \Phi, \ell \right)
      &= \te_{n+1} \etab^{-3} \te_{n+1} \etab^{2} \te_{n+1} \etab \left( \Phi, \ell \right) \nonumber \\
&= \te_{n+1} \etab^{-3} \te_{n+1}   \left(\etab^{3} \Phi, \ell \right) \nonumber \\
      &= \te_{n+1} \etab^{-3}   \left(\rr r_{0} \etab^{3} \Phi, \ell+1 \right) \label{eq:chapulines2_a} \\
      &= \te_{n+1}    \left(\etab^{-3} \rr r_{0} \etab^{3} \Phi, \ell+1 \right) \nonumber \\
&=     \left( \etab^{-3} \rr r_{0} \etab^{3} \Phi, \ell+1 \right) \label{eq:chapulines2_b}.
    \end{align}
  \end{subequations}
  Equation\nobreakspace \textup {(\ref {eq:chapulines2_a})} holds because $F_{0} \in \Psi = \kng^{3}\Phi$.
  Observe that $\kng^{-3}(\rr r_0 \Phi,\ell)$ is the root of its facet and clearly $F_{0} \not\in \etab^{-3}\rr r_0 \etab^{3} \Phi$, because $F_{0} \in \rr r_{0} \etab^{3}\Phi$.
  Moreover, $F_{1}$ cannot be the facet of $\kng^{-3}\rr r_0 \etab^{3} \Phi$ (Item\nobreakspace \ref {item:NMA} of  Proposition\nobreakspace \ref {prop:etaIsGood}).
  Hence, \eqref{eq:chapulines2_b} holds and this completes the proof of Item\nobreakspace \ref {item:chapulines22}.

  An analogous computation to that in Equation\nobreakspace \textup {(\ref {eq:chapulines2})} can be used to prove Item\nobreakspace \ref {item:chapulines23} (keeping in mind to replace $\ell+1$ by $\ell-1$).
   In order to prove Equation\nobreakspace \textup {(\ref {eq:chapulines2_b})} in this case we need to show that $F_3$ is not the facet of $\kng^{-3}\rr r_0 \etab^{3} \Phi$. However, according to Item\nobreakspace \ref {item:esAntipoda} of Proposition\nobreakspace \ref {prop:etaIsGood} this is only possible if $\rr r_0 \etab^{3}\Phi = r_{0}\hat{\Phi_{0}}$, and this, in turn, implies that $\Phi = r_{0}\hat{\Phi_{3}}$, which was excluded by hypothesis.

Let $\Phi' = \etab^{-3} \rr r_{0} \etab^{3} \Phi$. We will show that if $\ell \equiv 2 \pmod{3}$, then $(\Phi',\ell+1)$ satisfies the hypotheses of of Lemma\nobreakspace \ref {lem:chapulines2} and thus, Item\nobreakspace \ref {item:chapulines13} holds for $\Phi',\ell)$. First, note that the same argument used in the proof of Lemma\nobreakspace \ref {lem:chapulines1} can be used to prove that $\Phi' \neq r_{0} \hat{\Phi_{3}}$.
Further, $\ell+1 \equiv 0 \pmod{3}$.
It only remains to show that if $(\Psi',\ell') = \etab^{2} \te_{n+1} \etab (\Phi', \ell+1)$, then $F_{0} \in \Psi'$ and $\ell' \equiv 0 \pmod{3}$.
Just consider
\begin{subequations}
\begin{align}
	(\Psi',\ell')
	&= \etab^{2} \te_{n+1} \etab \left( \etab^{-3} \rr r_{0} \etab^{3} \Phi, \ell+1 \right) \nonumber \\
	&= \etab^{2} \te_{n+1} \left( \etab^{-2} \rr r_{0} \etab^{3} \Phi, \ell+1 \right) \nonumber \\
	&= \etab^{2} \left( \etab^{-2} \rr r_{0} \etab^{3} \Phi, \ell+1 \right) \label{eq:chapulines2_c} \\
	&= \left( \rr r_{0} \etab^{3} \Phi, \ell+1 \right). \nonumber
\end{align}
\end{subequations}
Recall that $F_{0} \in \etab^{3}\Phi$, which implies that $F_{0} \in  \rr r_{0} \etab^{3} \Phi$.
It follows that $\etab^{-2}\left( \rr r_{0} \etab^{3} \Phi, \ell+1 \right)$ is the root of its facet and this facet cannot be $F_{0}$ or $F_{1}$ (or $F_{3}$).
Therefore $\etab^{-2}\left( \rr r_{0} \etab^{3} \Phi, \ell+1 \right)$ is fixed by $\te_{n+1}$ which proves \eqref{eq:chapulines2_c}.
We just proved that $\ell' = \ell +1 \equiv 0 \pmod{3}$ and, as previously observed, $F_{0} \in  \rr r_{0} \etab^{3} \Phi = \Psi'$.

As before, a similar argument can be used to show that if $\ell \equiv 0 \pmod{3}$, then $(\Phi',\ell-1)$ satisfies the hypotheses of Item\nobreakspace \ref {item:chapulines12}.

Finally Item\nobreakspace \ref {item:chapulines24} follows from observing that
  \[\begin{aligned}
  \etab^{-3} \rr[\ell+1] r_{0} \etab^{3}
  (\etab^{-3} \rr r_{0} \etab^{3} \Phi) &= \Phi &\text{and} \\
  \etab^{-3} \rr[\ell-1] r_{0} \etab^{3} (\etab^{-3} \rr r_{0} \etab^{3} \Phi) &= \Phi,
    \end{aligned}
  \]
where in both cases we have used that $\rr[\ell-1]$, $\rr$ and $\rr[\ell+1]$ act the same way on flags containing $F_0$.
\end{proof}

\begin{lemma}\label{lem:chapulines3}
Let $(\Phi,\ell)$ be a flag, let $(\Psi,\ell') =  \etab^{-3} \te_{n+1} \etab^{2} \te_{n+1} \etab (\Phi,\ell)$. Suppose $F_0 \in \Psi$ and $\ell' \equiv i \pmod 3$ for some integer $i\in \{0,2\}$. Then the following hold:
\begin{enumerate}
    \item \label{item:chapulines31} $(\Phi,\ell) = (\Psi, \ell')$;
    \item \label{item:chapulines32} if $i = 2$ then $\mub(\Phi,\ell) = (\rr r_0 \Phi,\ell+1)$;
    \item \label{item:chapulines33} if $i = 0$ then $\mub(\Phi,\ell) = (\rr r_0 \Phi,\ell-1)$;
    \item \label{item:chapulines34}$|\left\langle \mub \right\rangle (\Phi,\ell)| = 2$.
\end{enumerate}

\end{lemma}

\begin{proof}
To prove Item\nobreakspace \ref {item:chapulines31} consider the flag
\[
    \etab^3(\Psi,\ell') = \te_{n+1} \etab^2 \te_{n+1} \etab (\Phi,\ell).
\]
Observe that $\etab^3(\Psi,\ell')$ is the root of its facet and that this facet is different from $F_0$ and $F_3$ (by Item\nobreakspace \ref {item:esVecina} of Proposition\nobreakspace \ref {prop:etaIsGood}).
Thus $\te_{n+1}$ acts trivially on $\etab^3(\Psi,\ell')$. It follows that
\[\begin{aligned}
\etab^{3}(\Psi,\ell')
&=\te_{n+1} \etab^3(\Psi,\ell') \\
&= \te_{n+1} \te_{n+1} \etab^2 \te_{n+1} \etab (\Phi,\ell) \\
&= \etab^2 \te_{n+1} \etab (\Phi,\ell),
\end{aligned}\]
thus
\[
(\Psi,\ell')
= \etab^{-3} \etab^{2} \te_{n+1} \etab (\Phi,\ell)
= \etab^{-1} \te_{n+1} \etab (\Phi,\ell).
\]

With a similar argument but now using that $\etab(\Psi,\ell')$ is the root of its facet one can prove that
\begin{align}
    \etab(\Psi,\ell') = \te_{n+1} \etab (\Phi,\ell) = \etab (\Phi,\ell),
\end{align}
which clearly implies that $(\Psi,\ell') = \etab^{-1}\etab^1 (\Phi,\ell) = (\Phi,\ell)$. 
This proves Item\nobreakspace \ref {item:chapulines31}.

In particular, $\mub(\Phi,\ell) = \te_{n+1}(\Phi,\ell)$ and by definition $\te_{n+1}(\Phi,\ell) = (\rr r_0 \Phi,\ell+j)$, where $j=1$ if $i=2$ and $j=-1$ if $i=0$.
Items\nobreakspace \ref {item:chapulines32} and\nobreakspace  \ref {item:chapulines33} follow.

Now we will prove Item\nobreakspace \ref {item:chapulines34}.
As before, we will just prove the case when $i=2$.
The case when $i=0$ follows from an analogous argument. Define $\Phi' = \rr r_{0} \Phi$ and $(\Psi',\ell') =\etab^{-3} \te_{n+1} \etab^{2} \te_{n+1} \etab(\Phi',\ell+1)$.
Notice that $F_{0} \in \Phi'$, which implies that $(\kng \Phi', \ell+1)$ and $(\kng^{3} \Phi', \ell+1)$ are the root of its respective facet and this facet is neither $F_{0}$, $F_{1}$ or $F_{3}$; therefore, they are fixed by $\te_{n+1}$.
Consider the following computation.
\[
\begin{aligned}
      \left( \Psi', \ell' \right)
      &= \etab^{-3} \te_{n+1} \etab^{2} \te_{n+1} \left(\etab \Phi', \ell+1 \right) \\
      &= \etab^{-3} \te_{n+1} \etab^{2}  \left( \etab \Phi', \ell +1 \right) \\
      &=  \etab^{-3} \te_{n+1}   \left(\etab^{3} \Phi', \ell+1 \right) \\
      &=  \etab^{-3}   \left( \etab^{3} \Phi', \ell+1 \right) \\
      &=(\Phi',\ell+1).
    \end{aligned}
\]
Now we clearly have that $F_{0} \in \Psi' = \Phi'$ and $\ell' = \ell+1\equiv 0 \pmod{3}$, as desired.
Finally the proof of Item\nobreakspace \ref {item:chapulines34} is complete by observing that
  \[\begin{aligned}
   \rr[\ell+1] r_{0}
  ( \rr r_{0}  \Phi) &= \Phi &\text{and} \\
   \rr[\ell-1] r_{0}  ( \rr r_{0} \Phi) &= \Phi,
    \end{aligned}
  \]
where as before, we used that $\rr[\ell-1]$, $\rr$ and $\rr[\ell+1]$ act the same way on flags containing $F_0$.
\end{proof} 

Lemmas\nobreakspace  \ref {lem:chapulines0} to\nobreakspace  \ref {lem:chapulines3}  prove that Lemma\nobreakspace \ref {lem:onlyThreeCopies} holds, which in turns proves Proposition\nobreakspace \ref {prop:noGroupAut}.

 \subsection{The intersection property} \label{sec:IP}
Now we are going to prove that the group $\groupExt$  satisfies the intersection property. 
The key to this proof is to use Lemma\nobreakspace \ref {lem:Lemma10}. 
In order to do so, we shall consider the following group elements
\begin{equation}\label{eq:gensEs}
\begin{aligned}
  \es_{1} &= \te_{1}, &&\\
  \es_{i} &= \te_{i-1}^{-1} \te_{i} && \text{if } 2 \leq i \leq n+1.
\end{aligned}
\end{equation}

Now for each $i$, $\es_{i}$ plays the role of $\sigma_{i}$ in Lemma\nobreakspace \ref {lem:Lemma10}.
Observe that
\begin{equation}
\langle \es_{1}, \dots, \es_{n} \rangle = \langle \te_{1}, \dots, \te_{n} \rangle \cong \monw(\cT) \cong \autp(\cT).
\end{equation}
Moreover, the elements $\es_{1}, \dots, \es_{n}$ act on $(\Phi, \ell)$ by \[\es_{i}(\Phi,\ell) = (s_{i} \Phi, \ell),\] where $s_{i} = r_{i-1}r_{i}$. 
Therefore, in order to use Lemma\nobreakspace \ref {lem:Lemma10} to prove the that $\groupExt$ satisfies the intersection property, we just need to show that
\begin{equation}
  \langle \es_{1}, \dots, \es_{n} \rangle \cap \langle \es_{j}, \dots, \es_{n+1} \rangle = \langle \es_{j}, \dots, \es_{n}\rangle
\end{equation}
for every $j$ such that $2 \leq j \leq n+1$ (note the change of rank).

The strategy to follow is to consider a certain flag $\Phi$ of  $\cT$ such that the intersection between the orbit of $(\Phi,1)$ under $\langle \es_{1}, \dots, \es_{n} \rangle $ and the orbit of $(\Phi,1)$ under $\langle \es_{j}, \dots, \es_{n+1} \rangle$ is precisely the orbit of $(\Phi,1)$ under  $\langle \es_{j}, \dots, \es_{n} \rangle$. The intersection property will follow from the fact that the action of $\langle \es_{1}, \dots, \es_{n} \rangle$ on the set $\Fw(\cT) \times \{1\}$ is free. The following results will focus on guaranteeing  the existence of such a flag, for which we will introduce some terminology.

Consider the base vertex $V_0$ of $\cT$, and let $\gamma_{1} \in \tras(\cT)$ be the translation by the vector $(1,0^{n-1})$. 
We define the \emph{line} $L$ as the orbit of the base vertex $V_0$ under $\langle \gamma \rangle$. 
We say that $\Psi=  (\tau, \vy, u)$ is \emph{perpendicular} to $L$ if $n \tau = 1$.
In geometric terms a flag is perpendicular to $L$ precisely when its $(n-1)$-face belongs to a hyperplane whose normal vector points in the same direction as $L$.

We will now prove some technical structural results about $\cT$.

\begin{rem}\label{rem:noBrincos}
If $\Phi$ is a flag whose vertex is in $L$, then the facet of $\Phi$ is neither $F_1$ nor $F_3$, nor does it belong to the set $H$ defined in Section\nobreakspace \ref {sec:chirExtRegToroids}. 
\end{rem}

\begin{lem}\label{lem:perpUnderLastS} 
  Let $\Phi$ be a flag of $\cT$ such that the vertex of $\Phi$ is in $L$ and assume that $\Phi$ is perpendicular to $L$. 
  Let $\gamma_1$ be the translation by the vector $e_1$ and let $\es_{n+1}$ be as in Equation\nobreakspace \textup {(\ref {eq:gensEs})}. 
  Then $\es_{n+1}(\Phi, 1) \in (\Phi, 1)\left\langle \gamma_{1} \right\rangle $. 
  Moreover,  all the flags in $\langle \es_{n+1} \rangle (\Phi, 1)$ are perpendicular to $L$ and their vertex is in $L$. 
\end{lem}
\begin{proof}
  First observe that the facet of $\Phi$ is not $F_1$ or $F_3$ and it does not belong to $H$ (see Remak\nobreakspace \ref {rem:noBrincos}), then \[\es_{n+1}(\Phi, 1) = \te^{-1}_{n}\te_{n+1} (\Phi,1) = \te^{-1}_{n}(\rr[1]r_0\Phi,1) = (r_{n} r_0 \rr[1]r_0 \Phi,1) = (r_{n} \rr[1] \Phi,1). \]
   Let $\Phi= (\sigma, \vx, t)$. Then we have
  \[\begin{aligned}
	  \es_{n+1}(\Phi, 1) &=   (r_{n} \rr[1] (\sigma, \vx, t),1)\\
	  &= (r_{n} (\sigma, \vx^{(1)}, t),1)\\ 
	  &= ((\sigma, \vx, t+x_{1}e_{1}),1),
  \end{aligned}\]
where the second equality holds because the facet of $\Phi$ is not in $H$, while the last equality holds because $n\sigma = 1$, since $\Phi$ is perpendicular to ${L}$.
  Finally observe that the flag $\Psi=(\sigma, \vx, t+x_{1}e_{1})$ also satisfies that its vertex belong to $L$ and it is perpendicular to $L$, hence the result follows.
\end{proof}

\begin{lem}\label{lem:translatedOfPerps} 
Let $\Phi= (\sigma,\vx,t)$ be a flag of $\cT$ such that the vertex  of $\Phi$ belongs to $L$ .  
  Let $j \geq 2$, $v \in \langle \es_{j}, \dots, \es_{n} \rangle$,
  and $(\Psi, \ell) = \es_{n+1} v \es_{n+1}^{-1}(\Phi, 1)$. 
  Then $\ell = 1$ and $\Psi$ is the image under a translation in direction $1$ of a flag in $\langle s_{j}, \dots, s_{n} \rangle \Phi$.
\end{lem}

\begin{proof}
  Let $v \in \langle \es_{j}, \dots, \es_{n}\rangle$. 
  Observe that $v$ fixes the second coordinate $(\Phi, 1)$ and thanks to Remak\nobreakspace \ref {rem:noBrincos}, so does $\es_{n+1}$.
  This indeed proves that $\ell = 1$.
  Moreover, we can think of $v$  as an element of $\langle s_{j}, \dots, s_{n} \rangle$.  
  Observe that we can write $v$ as 
  \[v=u_{r}s_{n} u_{r-1} s_{n} \cdots s_{n} u_{0},\] 
  with $u_{i} \in \langle s_{j}, \dots, s_{n-1} \rangle$ for $i \in \{0, \dots, r\}$. 
  We shall prove our result by induction over $r$.
  Assume that $r=0$, this implies that  \[\Psi = r_n \rr[1] u_0 \rr[1] r_n \Phi\] and recall that $u_0 \in \langle s_j, \dots, s_{n-1} \rangle$ commutes with $\rr[1]$, thus $\Psi = r_n u_0 r_n \Phi$ and $r_n u_0 r_n \in \langle s_j, \dots, s_{n} \rangle$. 

  Before proving the general case, let us prove the following fact: if $\Psi_1$ is a translate in direction $1$ of a flag in  $\langle s_{j}, \dots, s_{n} \rangle \Phi$, then so is  the flag $\Psi_2 = r_{n} \rr[1] s_{n} \rr[1] r_{n} \Psi_1$.
  Indeed, take $\Psi_1 = (\sigma, \vx, t)$ as above. 
  Notice that if $n \sigma \neq 1$, then $\rr[1] r_{n} \Psi_1 = r_{n} \rr[1] \Psi_1$. 
  In this situation \[r_{n} \rr[1] r_{n-1} r_{n} \rr[1] r_{n} \Psi_1 = r_{n} \rr[1] r_{n-1} r_{n}  r_{n} \rr[1] \Psi_1 = r_{n} r_{n-1}\Psi_1 = s_{n}^{-1} \Psi_1,\] where the second equality follows from the fact that $\rr[1]$ and $r_{n-1}$ commute (see Remak\nobreakspace \ref {rem:rchiquitas}). 
  Assume that  $n \sigma = 1$ and consider then the following computation:
	\[\begin{aligned}
		r_{n} \rr[1] r_{n-1} r_{n} \rr[1] r_{n} (\sigma, \vx, t) &= r_{n} r_{n-1} \rr[1]  r_{n} \rr[1] r_{n} (\sigma, \vx, t) \\ &=  r_{n} r_{n-1} \rr[1]  r_{n} \rr[1] (\sigma, \vx^{(1)}, t-x_{1}e_{1}) \\
		&= r_{n} r_{n-1} \rr[1]  r_{n} (\sigma, \vx, t-x_{1}e_{1})\\
		&= r_{n} r_{n-1} \rr[1]  (\sigma, \vx^{(1)}, t-x_{1}e_{1}- x_{1} e_{1}) \\
		&= r_{n} r_{n-1} (\sigma, \vx, t-x_{1}e_{1}- x_{1} e_{1}). \\
	 \end{aligned},\]
  where the first equality holds because $\rr[1]$ and $r_{n-1}$ commute (see Remak\nobreakspace \ref {rem:rchiquitas}).
  
  Let us come back to the general case. Assume that the result holds for any word  $u'_{r-1}s_{n} u'_{r-1} s_{n} \cdots s_{n} u'_{0}$ with $u_{i} \in \langle s_{j}, \dots, s_{n-1} \rangle$ and consider now $v$ as above.
  Define the flags $\Psi_1 = r_{n} \rr[1] u_{0} \rr[1] r_{n} \Phi = r_n u_0 r_n \Phi$ and  $\Psi_2 = r_{n} \rr[1] s_{n} \rr[1] r_{n} \Psi_1$. 
  From our previous discussion, $\Psi_2$ is the translate in the direction of $e_1$ of a flag in $\langle s_{j}, \dots, s_{n} \rangle \Phi$, in particular its vertex belongs to $L$.
  Notice that 
  \begin{align*}
    \Psi &= r_{n} \rr[1] v \rr[1] r_{n} \Phi \\
    &= (r_{n} \rr[1] u_{r} s_n \cdots s_n u_1 \rr[1] r_n)(r_{n} \rr[1] s_{n} \rr[1] r_{n})  (r_{n} \rr[1] u_{0} \rr[1] r_{n} \Phi) \\
    &= (r_{n} \rr[1] u_{r} s_n \cdots s_n u_1 \rr[1] r_n)(r_{n} \rr[1] s_{n} \rr[1] r_{n} \Psi_1) \\
    &= (r_{n} \rr[1] u_{r} s_n \cdots s_n u_1 \rr[1] r_n) \Psi_2. \\
  \end{align*} 
The result follows from our inductive assumption on $u_r s_n \cdots s_n u_1$.
\end{proof}

Using the previous results it is easy to obtain a description of the orbit of certain flags under the action of $\langle \es_{j}, \dots, \es_{n+1} \rangle$ for $j \geq 2$.

\begin{lem}\label{lem:IP_orbitGreatGens} 
  Let $\Phi$ be a flag of $\cT$ such that the vertex of $\Phi$ is in $L$. Assume that $\Phi$ is perpendicular to $L$. Let $j \geq 2$. Then
  \begin{equation}\label{eq:IP_orbitGreatGens}
	  \begin{split}
		  \langle \es_{j}, \dots, \es_{n+1}\rangle(\Phi,1) =\bigcup_{i=0}^{a-1} \langle \es_{j}, \dots, \es_{n}\rangle (\Phi,1)\gamma_{1}^{i}
	  \end{split}
  \end{equation}where $\gamma_{1} \in \tras(\cT)$ is the translation by the vector $e_{1}$. 
  In other words, the orbit of $(\Phi,1)$ under $\langle \es_{j}, \dots, \es_{n+1}\rangle$ is the union of translates in the direction of $e_1$  of the orbit of $\Phi$ under $\langle s_{j}, \dots s_{n}\rangle$. 
\end{lem}

\begin{proof}
   Observe that Lemma\nobreakspace \ref {lem:perpUnderLastS} proves that the right side of Equation\nobreakspace \textup {(\ref {eq:IP_orbitGreatGens})} is contained in the orbit of $(\Phi,1)$ under $\langle \es_{j}, \dots, \es_{n+1}\rangle$.

  To prove the other inclusion, it suffices to show that given a permutation $w \in \langle \es_j, \ldots, \es_{n+1} \rangle$ there exists $v \in \langle \es_j, \ldots, \es_{n} \rangle$ and $d \in \bZ$ such that $w(\Phi,1) = v(\Phi,1)\gamma_1^d$. First, observe that for some integer $r$ we can rewrite $w$ as $w = v_r \es_{n+1} v_{r-1} \es_{n+1} \ldots \es_{n+1}v_0$ where $v_i \in \langle \es_j, \ldots , \es_{n} \rangle$. We proceed by induction on $r$.

If $r=0$ then $w = v_0 \in \langle \es_j \ldots \es_{n} \rangle$, and we are done.

Let us now show that the case $r=1$ holds. Suppose that $r=1$ so that $w = v_1 \es_{n+1} v_0$ with $v_0,v_1 \in \langle \es_j , \ldots , \es_{n} \rangle$. Let $(\Psi,\ell) = \es_{n+1}(\Phi,1)$ and note that by Lemma\nobreakspace \ref {lem:perpUnderLastS} 
we have $\Psi = \Phi\gamma_1^{\epsilon}$ (for some $\epsilon \in \{\pm 1\}$) and $\Psi$ is perpendicular to $L$. Then
 \begin{align}
 \label{eq:orbitlemma1}
	  v_1 \es_{n+1} v_0 (\Phi, 1) =   v_1 \es_{n+1} v_0 \es_{n+1}^{-1}(\Psi, \ell) = v_1v_0'(\Psi,\ell)\gamma_1^{d} = v_1v_0'(\Phi,1)\gamma_1^{d+\epsilon}
  \end{align}

  for some $v_0' \in \langle \es_j, \ldots, \es_{n} \rangle$ and some integer $d$, where the second equality follows from Lemma\nobreakspace \ref {lem:translatedOfPerps}. Since $v_1v_0' \in \langle \es_j, \ldots, \es_{n+1} \rangle$, this concludes the case $r=1$.

  Now, suppose that the result holds for all $N < r$ and consider an element $w = v_r \es_{n+1} \ldots \es_{n+1}v_1 \es_{n+1}v_0$. By Equation\nobreakspace \textup {(\ref {eq:orbitlemma1})}, there exists $v_0' \in \langle \es_j,\ldots, \es_n \rangle$ and $d_1 \in \bZ$ such that
\begin{align}
   w(\Phi,1) 
   &= v_r \es_{n+1} \ldots \es_{n+1}v_1 \es_{n+1}v_0(\Phi,1) \nonumber \\
   &= (v_r \es_{n+1} \ldots \es_{n+1})v_1 \es_{n+1}v_0(\Phi,1) \nonumber \\
   &= (v_r \es_{n+1} \ldots \es_{n+1})v_1v_0'(\Phi,1)\gamma_1^{d_1} \nonumber \\
   &= (v_r \es_{n+1} \ldots \es_{n+1}v_1v_0')(\Phi,1)\gamma_1^{d_1} \label{eq:ladearriba}
  \end{align}

  Now let $w'= (v_r s_{n+1} \ldots s_{n+1}v_1v_0')$ so that Equation\nobreakspace \textup {(\ref {eq:ladearriba})} can be rewritten as $ w(\Phi,1) = w'(\Phi,1)\gamma_1^{d_1}$. Observe that $(\Phi,1)\gamma_1^{d_1}$ is perpendicular to $L$. Further, $w' \in \langle \es_j,\ldots,\es_{n+1} \rangle$ and $w'$ has less than $r$ factors equal to $\es_{n+1}$. Thus, by inductive hypothesis we have $w'(\Phi,1)\gamma_1^{d_1} = v(\Phi,1)\gamma_1^{d_1+d_2}$ for some $v \in \langle \es_j \ldots \es_{n} \rangle$ and $d_2 \in \bZ$.

\end{proof}

\begin{figure}\centering
\def\svgwidth{0.6\textwidth}
\begingroup \makeatletter \providecommand\color[2][]{\errmessage{(Inkscape) Color is used for the text in Inkscape, but the package 'color.sty' is not loaded}\renewcommand\color[2][]{}}\providecommand\transparent[1]{\errmessage{(Inkscape) Transparency is used (non-zero) for the text in Inkscape, but the package 'transparent.sty' is not loaded}\renewcommand\transparent[1]{}}\providecommand\rotatebox[2]{#2}\newcommand*\fsize{\dimexpr\f@size pt\relax}\newcommand*\lineheight[1]{\fontsize{\fsize}{#1\fsize}\selectfont}\ifx\svgwidth\undefined \setlength{\unitlength}{592.90782057bp}\ifx\svgscale\undefined \relax \else \setlength{\unitlength}{\unitlength * \real{\svgscale}}\fi \else \setlength{\unitlength}{\svgwidth}\fi \global\let\svgwidth\undefined \global\let\svgscale\undefined \makeatother \begin{picture}(1,1)\lineheight{1}\setlength\tabcolsep{0pt}\put(0,0){\includegraphics[width=\unitlength,page=1]{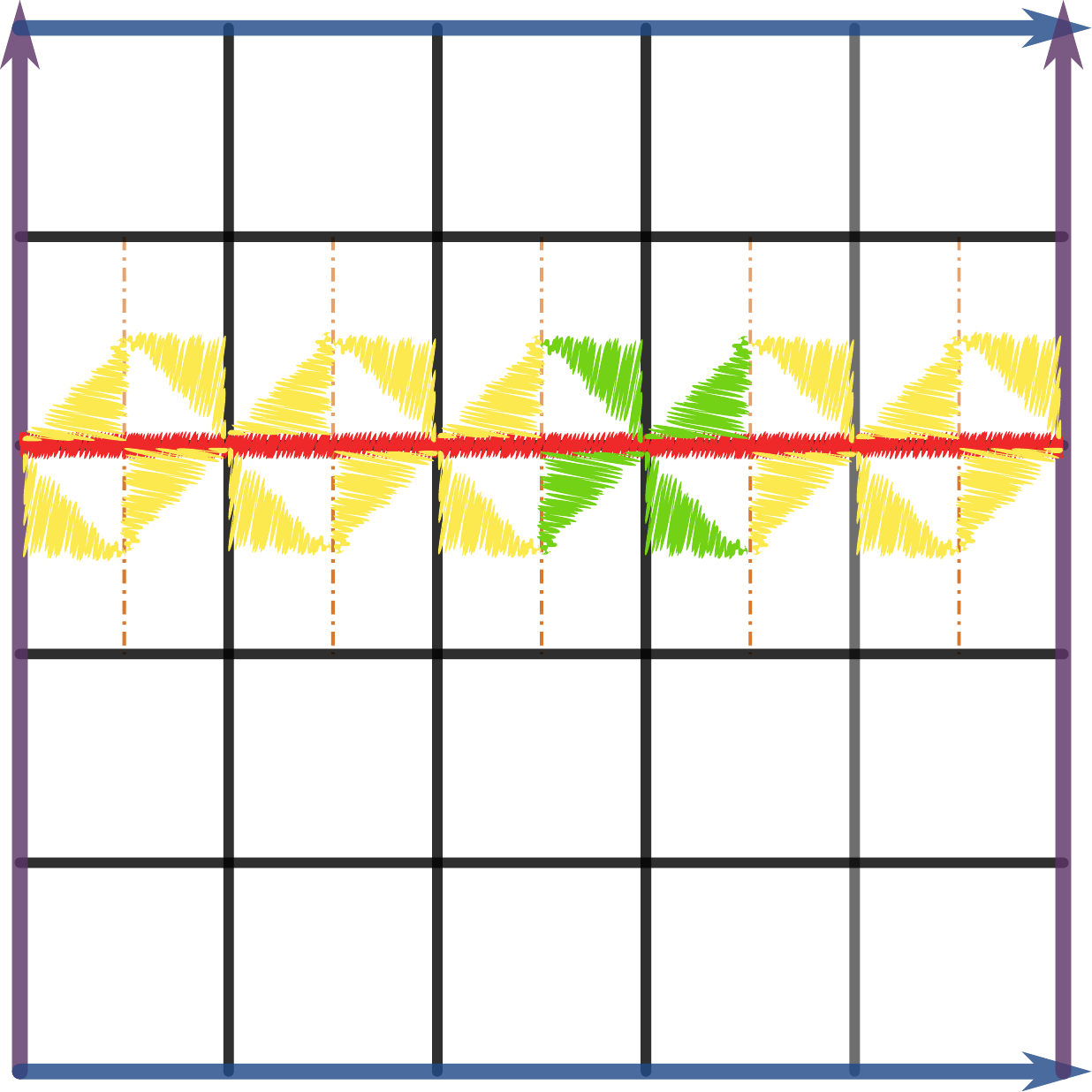}}\put(0.62789324,0.51107142){\color[rgb]{0,0,0}\makebox(0,0)[t]{\lineheight{1.25}\smash{\begin{tabular}[t]{c}$\Phi$\end{tabular}}}}\end{picture}\endgroup  \caption{The first coordinates of a pair $(\Phi,\ell)$ of Lemma\nobreakspace \ref {lem:IP_orbitGreatGens}.}
\end{figure}

\begin{lem}\label{lem:IP_FixedVertex} 
  Let $\Phi$ be a flag of $\cT$ perpendicular to $L$. Assume that the vertex of $\Phi$ belongs to $L$. Let $2\leq j \leq n$ and $w \in \langle \es_{1}, \dots, \es_{n} \rangle \cap \langle \es_{j}, \dots \es_{n+1} \rangle$. If $(\Phi', \ell) = w (\Phi, 1)$, then $\ell = 1$ and the vertex of $\Phi'$ is the same as the vertex of $\Phi$.
\end{lem}

\begin{proof}
We will slightly abuse language and say that $w$ fixes a vertex $V$, if for every $(\Psi,\ell)$ such that $\Psi$ has vertex $V$, the first coordinate of $w(\Psi,\ell)$ is also a flag with vertex $V$.

First, observe that since $w \in \langle \es_{1}, \dots, \es_{n} \rangle$, $w$ acts on the first coordinate of every pair $(\Psi,1)$ as a permutation $s \in \langle s_{1}, \dots, s_{n} \rangle \leq \monp(\cT)$. Furthermore, since $\cT$ is regular (and $s$ is in the monodromy group), we see that if $w$ fixes one vertex of $\cT$, then it must fix all vertices. 

Now, $w \in \langle \es_{j}, \dots \es_{n+1} \rangle$ and clearly all $s_{j}$ with $2 \leq j \leq n$ fix all vertices of $\cT$. This means that if $\es_{n+1}$ is not a factor of $w$, then $w$ fixes all vertices of $\cT$. Therefore, it suffices that we show that $\es_{n+1}$ fixes some vertex. In other words, we will show that for some flag $\Psi$ of $\cT$, the vertex of $(\Psi,\ell)$ is the same as the vertex of $\es_{n+1}(\Psi,\ell)$.

Let $X_H$ be the vertex of $\cT$ defined in  Definition\nobreakspace \ref {defn:hoyo} and let $\Psi_H$ be a flag with vertex $X_H$, so that the facet of $\Psi_H$ belongs to $H$. Observe that, by the definition of $\rho^\ell$ (Equation\nobreakspace \textup {(\ref {eq:rr})}), the vertex of $\rho^1 \Psi$ is $X_H$ as well. Now, consider the following computation.

\begin{align*}
    \es_{n+1} (\Psi_H,1) 
    &= \te_{n}^{-1}\te_{n+1}(\Psi_H,1) \\
    &= (r_nr_0\rho^{1}r_0 \Psi_H,1)\\
    &= (r_n\rho^{1}\Psi_H,1).
\end{align*}

Clearly, $r_n\rho^1\Psi$ has the same vertex as $\Psi$, as neither $r_n$ nor $\rho^1$ move $X_H$. That is, the vertex of $\es_{n+1} (\Psi,1)$ is $X_H$. Therefore, $w$ fixes all vertices of $\cT$. In particular, the vertex of $\Phi$ is the same as the vertex of $\Phi'$.
\end{proof}

Now we are ready to prove the key results that will help us to prove that the group $\gamma =\langle \es_{1}, \dots, \es_{n+1} \rangle$ satisfies the intersection property in Equation\nobreakspace \textup {(\ref {eq:intPropertyChiral})}.

\begin{lem}\label{lem:IP_orbitsSmallj}
  Let $\Phi$ be a flag of $\cT$ that is perpendicular to $L$ and whose vertex belongs to $L$. Let $2\leq j$, then
  \begin{equation}\label{eq:IP_orbits}
    \big(\langle \es_{1}, \dots, \es_{n} \rangle  \cap \langle \es_{j}, \dots, \es_{n+1} \rangle \big) (\Phi,1) = \langle \es_{j}, \dots, \es_{n} \rangle(\Phi, 1).
  \end{equation}
\end{lem}
\begin{proof}
  It is only necessary to prove that
  \[\big(\langle \es_{1}, \dots, \es_{n} \rangle \cap \langle \es_{j}, \dots, \es_{n+1} \rangle \big) (\Phi,1) \subset \langle \es_{j}, \dots, \es_{n} \rangle(\Phi, 1).\]

  By Lemma\nobreakspace \ref {lem:IP_orbitGreatGens},
  \begin{align*}
  \big(\langle \es_{1}, \dots, \es_{n} \rangle \cap \langle \es_{j}, \dots, \es_{n+1} \rangle \big) (\Phi,1)
  &\subset \langle \es_{j}, \dots, \es_{n+1} \rangle (\Phi,1) \\
  &\subset \left(\bigcup_{i=0}^{a-1} \langle \es_{j}, \dots, \es_{n}\rangle (\Phi,1)\gamma_{1}^{i}\right),
  \end{align*} 
Now, by Lemma\nobreakspace \ref {lem:IP_FixedVertex}, if $(\Psi,1)$ is an element of $\langle \es_{j}, \dots, \es_{n+1} \rangle (\Phi,1)$, then $\Psi$ has the same vertex as $\Phi$. This implies that
   \[\langle \es_{1}, \dots, \es_{n} \rangle (\Phi,1) \cap \langle \es_{j}, \dots, \es_{n+1} \rangle (\Phi,1) \subset \langle \es_{j}, \dots, \es_{n}\rangle (\Phi,1). \qedhere
  \]
\end{proof}

As explained before Lemma\nobreakspace \ref {lem:IP_orbitsSmallj} offers the conditions to prove the intersection property for the group $\langle \es_{1}, \dots, \es_{n+1} \rangle$.

\begin{prop}\label{prop:IntersectionPropertyExtension} 
  Let $\es_{1}, \dots, \es_{n+1}$ be the group elements defined in Equation\nobreakspace \textup {(\ref {eq:gensEs})}. Let $j \geq 2$. Then
  \begin{equation}\label{eq:IP_groups}
    \langle \es_{1}, \dots, \es_{n} \rangle  \cap \langle \es_{j}, \dots, \es_{n+1} \rangle  = \langle \es_{j}, \dots, \es_{n} \rangle.
  \end{equation}
\end{prop}

\begin{proof}
  If $j = n+1$, then there is nothing to prove. If $j < n+1$, then let $\Phi$ be as in the hypothesis of Lemma\nobreakspace \ref {lem:IP_orbitsSmallj}. Let $w \in \langle \es_{1}, \dots, \es_{n} \rangle  \cap \langle \es_{j}, \dots, \es_{n+1} \rangle$ and observe that
\begin{align*}
    \big(\langle \es_{1}, \dots, \es_{n} \rangle  \cap \langle \es_{j}, \dots, \es_{n+1} \rangle \big) (\Phi,1) = \langle \es_{j}, \dots, \es_{n} \rangle(\Phi, 1)
  \end{align*}
  
 by Lemma\nobreakspace \ref {lem:IP_orbitsSmallj}. This implies that there exists $w' \in \langle \es_{j}, \dots, \es_{n} \rangle$ such that $w(\Phi,1) = w'(\Phi,1)$. Observe that the action of the group $\langle \es_{1}, \dots, \es_{n} \rangle$ on the set \[\{(\Psi,1) : \Psi\ \text{is a white flag of}\ \cT\}\] is equivalent to the action of $\monp(\cT)$ on the set of white flags of $\cT$. In particular, this action is free. Now, since both $w$ and $w'$ belong to the group $\langle \es_{1}, \dots, \es_{n} \rangle$ and $w(\Phi,1) = w'(\Phi,1)$, then $w=w'$, and we have that $w \in \langle \es_{j}, \dots, \es_{n}\rangle$.
  The other inclusion is obvious.
\end{proof}

As a consequence of  Propositions\nobreakspace \ref {prop:chiralExtRelations},  \ref {prop:noGroupAut} and\nobreakspace  \ref {prop:IntersectionPropertyExtension} we have the following.

\begin{thm}\label{thm:elChido}
Let $\cT$ be a regular toroid $\typetor_{\va}$ where $\va = (a^{k}, 0^{n-k})$ with $a \geq 6n +1$ and $k \in \{1,2,n\}$.
Let $\es_{1}, \dots, \es_{n+1}$ be the corresponding group elements defined in Equation\nobreakspace \textup {(\ref {eq:gensEs})}. Then the group $\groupExt = \langle \es_{1}, \dots, \es_{n+1} \rangle$ is the automorphism group of a chiral $(n+2)$-polytope $\cP$ whose facets are isomorphic $\cT$.
\end{thm}

The abstract polytopes constructed in Theorem\nobreakspace \ref {thm:elChido} prove Theorem\nobreakspace \ref {thm:elNoTanChido}.
 
\section*{Acknowledgements}
This research was partially developed when the first author was a Ph.D. student in Centro de Ciencias Matemáticas at UNAM Mexico. 

The second author gratefully acknowledges financial support from the F\'ed\'eration Wallonie-Bruxelles -- Actions de Recherche Concert\'ees (ARC Advanced grant).

Both authors are grateful to Daniel Pellicer for his insights and contributions at the beginning of this project.

\printbibliography

\end{document}